\documentclass[onefignum,onetabnum,final]{siamart220329}

\usepackage{comment}
\usepackage{lipsum}
\usepackage{amssymb}
\usepackage{tikz} 
\usepackage{amsfonts}
\usepackage{dsfont}
\usepackage{mathtools} 

\usepackage{graphicx}
\usepackage[caption=false]{subfig}
\usepackage{epstopdf}
\usepackage{algorithmic}
\usepackage{amsopn}
\usepackage{soul}          
\usepackage{enumitem}
\ifpdf
  \DeclareGraphicsExtensions{.eps,.pdf,.png,.jpg}
\else
  \DeclareGraphicsExtensions{.eps}
\fi

\definecolor {uclablue}   {RGB} {39,116,174}
\definecolor {masongreen} {RGB} {0,102,51}

\newcommand{\bigoh}[1]{\mathcal{O}\left( #1 \right)}

\newcommand{\etaeps}[0]{\eta_{\epsilon}}
\newcommand{\thetaeps}[0]{\theta_{\epsilon}}
\newcommand{\zeps}[0]{z_{\epsilon}}
\newcommand{\teps}[0]{t_{\epsilon}}

\newcommand{\MoscoConv}{\overset{M}{\longrightarrow}}

\newcommand{\GammaConvs}{\overset{\Gamma_s}{\longrightarrow}}

\DeclareMathOperator*{\argmin}{argmin}

\newcommand{\expval}[1]{\mathbb{E}\left[ #1 \right] }

\newcommand{\RR}[0]{\mathbb{R}}
\newcommand{\PP}[0]{\mathbb{P}}
\newcommand{\EE}[0]{\mathbb{E}}
\newcommand{\NN}[0]{\mathbb{N}}
\newcommand{\RRext}[0]{\overline{\mathbb{R}}}
\newcommand{\Rho}{\mathrm{P}}

\newcommand{\cA}[0]{\mathcal{A}}

\newcommand{\cE}[0]{\mathcal{E}}
\newcommand{\cV}[0]{\mathcal{V}}

\DeclareMathOperator*{\esssup}{ess\,sup}

\newsiamremark{remark}{Remark}
\newsiamremark{hypothesis}{Hypothesis}
\crefname{hypothesis}{Hypothesis}{Hypotheses}
\newsiamthm{claim}{Claim}
\newsiamthm{assumption}{Assumption}

    \usepackage[textsize=small]{todonotes}
\headers{Rockafellian relaxation for PDECO}{H.\ Antil, S.\ P. Carney, H.\ D\'{i}az, and J.\ O. Royset}

\title{Rockafellian relaxation for PDE-constrained optimization with distributional ambiguity
\thanks{
\funding{The first three authors are partially supported by NSF grant DMS-2110263, 
Air Force Office of Scientific Research (AFOSR) under Award NO: FA9550-22-1-0248, 
and Office of Naval Research (ONR) under Award NO: N00014-24-1-2147. The fourth author 
is supported by ONR under Award NO: N00014-24-1-2318.}
}}

\author{
Harbir Antil\footnotemark[1] \and
Sean P. Carney\footnotemark[2] \and 
Hugo D\'{i}az\footnotemark[3]  \and
Johannes O. Royset\footnotemark[4] 
} 

\ifpdf
\hypersetup{
  pdftitle={Rockafellian relaxation for PDECO with distributional ambiguity},
  pdfauthor={Harbir Antil, Sean P. Carney, Hugo D\'{i}az, and Johannes O. Royset}
}
\fi

\graphicspath{ {figs/}{./figs/} }
\DeclareGraphicsExtensions{.pdf}
\begin{document}

\maketitle
\renewcommand{\thefootnote}{\fnsymbol{footnote}}
\footnotetext[1]{Department of Mathematical Sciences and Center for Mathematics and Artificial Intelligence, 
George Mason University, Fairfax, VA, 22030, USA; hantil@gmu.edu}
\footnotetext[2]{Department of Mathematics,  Union College, Schenectady, NY 12308, USA; carneys@union.edu}
\footnotetext[3]{Department of Mathematics, North Carolina State University, Raleigh, NC, 27695, USA; hsdiazno@ncsu.edu}
\footnotetext[4]{Daniel J.\ Epstein Department of Industrial and Systems Engineering,  University of Southern California, Los Angeles, CA 90089, USA; royset@usc.edu}

\begin{abstract}
Stochastic optimization problems are generally known to be unstable
to the form of the underlying uncertainty. 
A framework is introduced
for optimal control problems with 
partial differential equations as constraints
that is robust to inaccuracies in the precise form of the 
problem uncertainty.
The framework is based on problem relaxation and involves optimizing 
a bivariate, ``Rockafellian'' objective functional that features both 
a standard control variable and an additional perturbation 
variable that handles the distributional ambiguity. 
In the presence of distributional corruption, 
the Rockafellian objective functionals are shown in 
the appropriate settings
to $\Gamma$-converge to 
uncorrupted objective functionals in the limit of vanishing
corruption. Numerical examples illustrate the framework's utility 
for outlier detection and removal and for variance reduction.

\end{abstract}



\section{Introduction}
\label{sec:intro}
Many stochastic optimization problems constrained by partial differential equations (PDEs)
take the form
\begin{equation}\label{eq:og_objective}
\min_{z \in Z_{\rm ad}} \varphi(z), \qquad 
\varphi(z) = f_0(z) + \expval{ g(s(\xi, z))} = f_0(z) + \int_{\Xi} g(s(\xi, z)) \, d\PP(\xi)
\end{equation}
where
the deterministic control variable $z$ belongs to some admissible set $Z_{\rm ad}$. 
Precise definitions are given below, but generally speaking, the function $s$ maps the control variable 
to the solution of the underlying PDE constraint, 
while $g$ maps the PDE solution to some quantity-of-interest, 
and $f_0(z)$ can be a control penalty or regularization term.
As the precise form of one or more of the PDE input data
may be unknown, the solution map $s$ is parameterized by some random quantity $\xi$ that belongs to 
a sample space $\Xi$ and follows a hypothesized distribution $\PP$.

Problem uncertainty can arise from imprecise knowledge of the constraining equation's
forcing terms, boundary and initial conditions, geometry, and coefficients within 
the differential operator, for example. In practice, the form of the uncertainty itself
is often ambiguous; i.e.\ the sample space $\Xi$ and probability measure $\PP$ are
themselves uncertain. Typically one makes an ansatz based on empirical knowledge, which 
may be unsettled due to measurement error or adversarial corruption, for example.

One standard approach to guard against such ``meta-uncertainty'' is to consider distributionally 
robust optimization (DRO) formulations, in which one minimizes an expected value over a 
collection $\{\PP_i\}_{i \in \mathcal{I}}$ of plausible measures
(called the ambiguity set), 
among which it is hoped that the correct, unknown distribution resides. 
In conservative approaches, 
one considers worst-case scenarios, which leads to the minimax problem 
\begin{equation}\label{eq:DRO_formulation}
\min_{z \in Z_{\rm ad}} \Big( f_0(z) + \sup_{i \in \mathcal{I}} \int_{\Xi} g(s(\xi, z)) \, d\PP_{i}(\xi) \Big).
\end{equation}
Such ``higher'', conservative formulations
are sensible for applications for which there exist outlier events
that have catastrophic, severe societal impacts and are to be guarded against at all costs. 
In practice, the full minimax problem \eqref{eq:DRO_formulation} may be intractable, and 
typically some kind of approximation is needed for computation, for example based on Taylor expansions of 
the map $g$ with respect to the random parameter $\xi$. The literature on DRO is vast; 
in the particular context
of PDE constrained
optimization (PDECO), some examples include  
\cite{lass2017model,kouri2017measure,kolvenbach2018approach,milz2022approximation} and the references therein. 
We also note that DRO is closely connected to the 
use of risk measures \cite{royset2025risk}, for example the conditional value-at-risk
\cite{kouri2016risk}.

In other 
applications, outlier events might not be catastrophic, and hence 
the worst-case, DRO approach \eqref{eq:DRO_formulation} may be too conservative.
 Especially for scenarios with tight performance requirements, a more optimistic 
approach is warranted in which best-case
scenarios are considered, which is termed distributionally favorable, or 
distributionally optimistic optimization (DOO). Optimistic approaches are generally 
based on problem relaxation, rather than restriction; for problems with distributional 
ambiguity, they allow one to keep in play decisions that might be optimal
under the true, uncorrupted distribution. In contrast, the DRO approach rules them out. 
A similar dichotomy between best-case and worst-case formulations appears in bilevel 
optimization; see e.g.\ \cite{DempeZemkoho.20}.

One strong reason to consider optimistic approaches is that, in general,  
optimization problems are prone to instabilities in the sense 
that small changes to problem data may cause large changes to their solutions, e.g.
\cite{royset2020stability,bonnans2013perturbation,levy2000stability}. 
This is true in particular for stochastic optimization, as discussed for example
in \cite{romisch1991distribution,romisch2007stability}. 
Here and throughout the manuscript, we informally define ``instability'' to 
mean that a small change in a probability distribution can cause a large 
change in an optimization problem's solution, i.e.\ 
the minimizer to a perturbed stochastic program can differ substantially from the 
minimizer to the unperturbed one. This can cause an upwards shift in the 
value of the objective function.  
An example of this instability for a simple one-dimensional 
problem is described below in \cref{subsec:motivation}, 
while numerous examples in the context of linear programming
can be found in \cite{ben2000robust}. 
See \cite{royset2024rockafellian} for additional 
examples in a variety of contexts.  
Multiple examples for optimal control problems  
of the form \eqref{eq:og_objective} are described in 
 \cref{subsec:numerical_examples}. 

Applied to problems with these kinds of instabilities, the DRO approach necessarily leads to 
higher objective function values because of the 
supremum in \eqref{eq:DRO_formulation}.
In contrast, DOO approaches seek lower, best-case scenarios, for example by 
replacing the supremum in \eqref{eq:DRO_formulation} with an infimum.
In general, this can lead to an intractable problem; for example, 
the authors in \cite{jiang2024distributionally} show that computing the inner infimum 
is NP-hard for a recourse function represented as a linear program with objective uncertainty. 
This is not always the case, however, and 
the DOO approach has been successfully used in the contexts of 
statistical learning 
\cite{norton2017optimistic,agarwal2020optimistic,song2020optimistic,cao2021contextual,gotoh2023data},
Bayesian optimization \cite{nguyen2019calculating,nguyen2019optimistic,nguyen2020robust}, 
and outlier analysis \cite{aravkin2020trimmed,zheng2021trimmed,narasimhan2023learning}. 
Besides connecting DOO to techniques from robust statistics and outlier analysis, the 
recent work in \cite{jiang2024distributionally} also integrates DOO and DRO together 
to derive out-of-sample performance guarantees. 
See also the recent works \cite{blanchet2024distributionally} 
and \cite{blanchet2024automatic} for more details on the 
connection between DRO and robust statistics and outlier detection and removal,
respectively. 

In the current work we present an optimistic formulation for PDECO 
under distributional uncertainty, which to the best of the authors' knowledge 
has not previously been considered. 
Rather than simply replacing the supremum in \eqref{eq:DRO_formulation} with an infimum 
(or a suitable reformulation thereof), the method seeks 
best-case scenarios based on problem perturbation, as recently proposed in 
\cite{royset2024rockafellian}.

In this approach, 
the original objective functional in \eqref{eq:og_objective} is generalized to a bivariate, 
``Rockafellian'' \cite{rockafellar1963convex,convex_analysis} objective functional
that depends on both the original control variable $z$ and additional
perturbation parameters. The Rockafellian is chosen so that it is 
equivalent to, or ``anchored at'', the 
original objective functional whenever the perturbation variable equals zero. 
The original stochastic optimization problem is thus embedded in a family of 
perturbed problems; the key point is that they 
are better conditioned to data corruptions and meta-uncertainty than both
\eqref{eq:og_objective} and DRO approaches.
For finite dimensional problems, 
the authors in 
\cite{royset2024rockafellian} show that  
Rockafellian objective functionals $\Gamma$-converge\footnote{
Note that in the literature, $\Gamma$-convergence is sometimes termed epi-convergence
\cite{feinberg2023epi,kouri2020epi}.}
\cite[Definition 12.1.1]{attouch2006variational}
to the original one whenever the size of 
the corruption on $\PP$ or $\Xi$ vanishes. As is well-known, $\Gamma$-convergence preserves 
sequences of converging minimizers.

There are two primary objectives for the present work. 
The first is to 
extend the mostly finite-dimensional theoretical 
results in \cite{royset2024rockafellian} to the infinite-dimensional setting of PDECO. 
Whenever the original objective functional \eqref{eq:og_objective} is corrupted, we 
prove 
$\Gamma$-convergence of suitably defined Rockafellians  
to \eqref{eq:og_objective} as the size 
the corruption vanishes. 
The $\Gamma$-convergence is shown with respect to the strong topology of 
a Banach space.
We show results both for the case of a corrupted probability 
density function, and for the case of a corruption to the support of a 
probability distribution.
We consider also a special example of the former case, namely when the probability 
space is both discrete and finite-dimensional, 
for example, when using a sample-average approximation (SAA) for computing 
expectation values.
In this setting we show Mosco convergence
whenever the corruption to the discrete probabilities vanish. 

A closely related work is from \cite{feinberg2023epi}, where the authors show 
 $\Gamma$-convergence of expectation functions under both varying measures and varying
integrands. In the context of PDECO (and using notation from the current work), 
the authors establish conditions under which functionals of the form 
$$
\EE_{\PP_{\epsilon}}\big[ g_{\epsilon}\big(s_{\epsilon}(\xi,z)\big)\big] 
= \int_{\Xi }g_{\epsilon}\big(s_{\epsilon}(\xi,z)\big) \, d\PP_{\epsilon}(\xi)
$$
$\Gamma$-converge to the expectation function in \eqref{eq:og_objective}. 
Here the $\PP_{\epsilon}$ are approximate probability measures that, in 
some sense, get close to $\PP$,  
 while $s_{\epsilon}$ and $g_{\epsilon}$ are approximations 
to $s$ and $g$ arising, for example, from numerical discretizations. 
Also related is the recent work in 
\cite{chen2023performance}, where optimality gaps for optimal controls in PDECO 
are derived 
under various kinds of general inaccuracies, for example, finite dimensional 
approximations, sample average approximations, or smooth approximations of nonsmooth 
functions. 

The results in \cite{feinberg2023epi} are, in some sense, 
more general than those contained here; 
$g$ and $s$ are maps between metric spaces, and $\Gamma$-convergence is shown whenever
the probability measures $\PP_{\epsilon}$ converge weakly to $\PP$ as $\epsilon \downarrow 0$. 
As mentioned above, this implies that accumulation points of minimizers of the approximating functionals
are minimizers of the $\Gamma$-limit. 
Although the metric space setting appears sufficiently abstract to handle the infinite-dimensional
setting of PDECO, often one is unable to show that there exist
accumulation points of minimizers in a metric topology for the function spaces (e.g.\ $L^p$, Sobolev)
that typically arise in PDECO.
An exception to this is when one can show compactness results for the solution operator 
of the PDE at hand. Recent efforts to ensure norm compactness of sets of critical points that 
appear in \cite{milz2024sample} may allow for wider application of the results from \cite{feinberg2023epi}, 
but such compactness tends to require additional regularity assumptions, as seen in \cite{milz2024sample}.  

The second objective of the present work is to 
showcase with numerical examples the practical utility of Rockafellian relaxation for 
PDECO problems; 
the method enables recovery of the optimal control to an \emph{uncorrupted} 
optimization problem even when solving with corrupted data. 
A simple demonstration for a one dimensional stochastic program can be found 
below in \cref{subsec:motivation}. 
The numerical examples 
in \cref{subsec:numerical_examples} illustrate the method's potential for 
both outlier detection (and subsequent
removal) and variance reduction in the presence of corruptions to probability 
densities and the support of a probability density, respectively. 

The rest of the paper is organized as follows. 
\Cref{sec:motiv_and_prelim} motivates
the Rockafellian relaxation framework 
with a simple example and introduces some preliminary definitions
and assumptions;
\cref{sec:general_theory} develops the general theory  
for PDECO problems of the form \eqref{eq:og_objective}. 
In \cref{sec:elliptic_pdes}, the conditions necessary for the general theory to hold 
are verified in the case of stochastic elliptic PDE constraints. 
\Cref{subsec:numerical_examples} then describes numerical examples of Rockafellian relaxation for PDECO, 
followed by 
a brief discussion of the important role played by the relaxation parameter inherent to the approach
in \cref{subsec:theta}.
After a discussion of the 
computational cost of the method in 
\cref{subsec:cost}, 
\Cref{sec:conclusions} summarizes the results and concludes with some possible 
directions for future research. 


\section{Motivation and preliminaries}
\label{sec:motiv_and_prelim}

   \subsection{Motivation}
   \label{subsec:motivation}
   As a simple example of a stochastic optimization problem that is unstable to 
perturbations in the underlying probability distribution, we borrow 
from \cite[Example 2.1]{royset2024rockafellian}.
Consider
the one-dimensional stochastic program $\min_{x \in [0,1]} \varphi(x)$, 
where
\begin{equation} \label{eq:1d_pure}
\varphi(x) = \expval{g(x,\xi)}, \qquad g(x,\xi) = (1-x)/2 + \xi x, 
\end{equation}
and
$\PP[\xi=0] = 1$. The
global minimizer is $x^{\star} = 1$. 

If instead, however, for some 
$0 < \epsilon  \ll 1$ 
we consider the corrupted random variable 
$\xi_{\epsilon}$ whose law is given by 
$p_{\epsilon,1}:= \PP[\xi_{\epsilon} = 0] =  1-\epsilon$ and 
$p_{\epsilon,2}:= \PP[\xi_{\epsilon} = 1/\epsilon] = \epsilon$, 
then
 the global minimizer of 
\begin{equation} \label{eq:1d_corrupted}
\varphi_{\epsilon}(x) = \expval{g(x,\xi_{\epsilon})} =  
\frac12 (1-x) + 0 \cdot p_{\epsilon,1} \, x  + 1/\epsilon \cdot p_{\epsilon,2} \, x
= \frac12(1+x) 
\end{equation}
on $[0,1]$ is 
 $x_{\epsilon}^{\star} = 0$. 
 Since $p_{\epsilon,1} \to 1$ and $p_{\epsilon,2} \to 0$ as $\epsilon \downarrow 0$, 
$\xi_{\epsilon}$ converges in distribution to $\xi$. However, no matter how small the perturbation
$\epsilon$ is, the minimizers $x_{\epsilon}^{\star}$ never approach $x^{\star}$.

To recover the uncorrupted minimizer using Rockafellian relaxation introduced 
in \cite{royset2024rockafellian}, we construct a bivariate function 
$\Phi_{\epsilon}: [0,1] \times \RR^2 \to \RRext:= \RR \cup \{-\infty,\infty\}$. Let $p_{\epsilon} = 
(p_{\epsilon,1}, p_{\epsilon,2})$, $t = (t_1,t_2)$, and let 
$$
\Delta := \Big\{ q \in \RR^2: \, \, q_1+q_2 = 1  \text{ and } 0 \le q_i \le 1  \text{ for } i \in \{1,2\} \Big\}
$$
denote the set of probability vectors on $\RR^2$. For 
some $\thetaeps > 0$, consider 
\begin{align}
\Phi_{\epsilon}(x,t) 
&= \frac12 (1-x) + 0 \cdot (p_{\epsilon,1} + t_1) x  +  1/\epsilon \cdot (p_{\epsilon,2} + t_2) x
+  \frac{\theta_{\epsilon}}{2} \| t\|_2^2 + \iota_{\Delta}(p_{\epsilon} + t) \nonumber \\
&= \frac12 (1-x) + (\epsilon + t_2) x/\epsilon  
+  \frac{\theta_{\epsilon}}{2} \| t\|_2^2 + \iota_{\Delta}(p_{\epsilon} + t) , \label{eq:Phieps_1d}
\end{align}
where $\| \cdot \|_2$ is the Euclidean norm and the indicator function 
\begin{equation}\label{eq:indic_defn}
\iota_{\Delta} (q) = \begin{cases} 0, \qquad &q \in \Delta \\ \infty, \qquad &q \notin \Delta. \end{cases}
\end{equation}
From elementary calculus, the quadratic program 
$$
\min_{x \in [0,1], t \in \RR^2} \Phi_{\epsilon}(x,t) 
$$
has a global minimum at $x_{\epsilon}^{\star} = 1$ and $t_{\epsilon}^{\star} = (\epsilon,-\epsilon)$, 
so long as $\thetaeps < (\epsilon/2)^{-2}$.

Thus, by solving the relaxed problem $\min_{x,t} \Phi_{\epsilon}(x,t)$ with 
additional perturbation variable $t$, 
the problematic data point $\xi_{\epsilon}=1/\epsilon$ that causes the 
upwards shift from $\min_x \varphi(x)$ to $\min_x \varphi_{\epsilon}(x)$  
is identified and removed. 
We observe a similar utility in 
analogous numerical examples in the context of stochastic PDECO, 
which we describe below in \cref{subsec:numerical_examples} after 
the 
theoretical developments in \cref{sec:general_theory}.


   \subsection{Preliminaries}
   \label{subsec:preliminaries}
   To develop
a general theory of Rockafellian relaxation for PDECO problems 
of the form \eqref{eq:og_objective} in the presence
of distributional ambiguity, we first introduce some
preliminary definitions and assumptions that will generally be used 
throughout \cref{sec:general_theory}. 

Let $(U,\|\cdot\|_U)$ and $(Z,\| \cdot \|_Z)$ be two Banach spaces, and let $(\Xi,\cA, \PP)$ be a probability space. 
For maps $f_0: Z \to \RRext$, $s: \Xi \times Z \to U$ and $g: U \to \RRext$
we invoke the following assumptions: 

\begin{assumption}[Properties of the solution map $s=s(\xi,z)$]\label{asmpt:1}
\begin{enumerate}
\item $s(\cdot, z): \Xi \to U$ is 
$\mathcal{A}$ measurable $\forall z \in Z$.

\vskip0.5em
\item If
$z_{\epsilon} \rightharpoonup z$ in $Z$ as $\epsilon \downarrow 0$, then 
$s(\xi, z_{\epsilon}) \rightharpoonup s(\xi, z)$ in $U$ 
a.s.\ in $\Xi$. 

\end{enumerate}
\end{assumption}

\begin{assumption}[Properties of $f_0$ and $g$]\label{asmpt:2}
\begin{enumerate}
\item $f_0$ is proper: $f_0(z) > -\infty$ $\forall z \in Z$ and 
$f_0(z) < \infty$ for some $z \in Z$. 

\vskip0.5em
\item 
Both $f_0$ and $g$ are weakly sequentially lower semi-continuous (lsc) maps: 
\begin{align*}
&z_{\epsilon} \overset{Z}{\rightharpoonup} z \implies \liminf_{\epsilon \downarrow 0} f_0(\zeps) \ge f_0(z) \\
&u_{\epsilon} \overset{U}{\rightharpoonup} u \implies \liminf_{\epsilon \downarrow 0} g(\zeps) \ge g(z). 
\end{align*}

\item There exists some 
$\gamma \in \RR$ 
such that $\forall u \in U$, 
$
g(u) \ge \gamma.
$

\end{enumerate}
\end{assumption}


We next briefly recall
the definitions of Mosco and $\Gamma$-convergence, as well as a standard result that follows from the definitions. 
Here we follow the definition of $\Gamma$-convergence given in
\cite[Definition 12.1.1]{attouch2006variational}, 
which involves maps from metric spaces to the extended real numbers. 
In this work we are primarily interested in Banach spaces, and in particular, 
we use the notion of $\Gamma$-convergence with respect to the strong (i.e.\ norm) topology.
For Mosco convergence, we follow \cite[Definition 17.4.6]{attouch2006variational}.

\begin{definition}[Mosco and $\Gamma_s$-convergence]\label{defn:mosco}
Let $(X, \|\cdot\|_X)$ be a Banach space, let $g: X \to \RRext$, and let 
$(g_{\epsilon})_{\epsilon\in\RR_+}$ be a sequence of maps from $X$ to $\RRext$ indexed by
$\epsilon >0$. 
The sequence $(g_{\epsilon})_{\epsilon\in\RR_+}$ 
Mosco converges to $g$,
$g_{\epsilon} \MoscoConv g$,
as $\epsilon \downarrow 0$ if
\begin{align*}
&\text{(i) } Consistency: \, \forall x \in X, \exists \, (x_{\epsilon})_{\epsilon\in\RR_+} \text{ such that } x_{\epsilon} \to x \text{ and } \limsup_{\epsilon \downarrow 0 } g_{\epsilon}(x_{\epsilon}) \le g(x). \\
&\text{(ii) Stability: for all sequences } x_{\epsilon} \rightharpoonup x \text{ we have: } \liminf_{\epsilon \downarrow 0 } g_{\epsilon}(x_{\epsilon})   \ge g(x).  
\end{align*}
The sequence $(g_{\epsilon})_{\epsilon\in\RR_+}$ $\Gamma$-converges to $g$ in 
the strong sense, 
denoted by $g_{\epsilon}\GammaConvs g$, if condition (i) holds and 
$$
\text{for all sequences } x_{\epsilon} \to x \text{ we have: } \liminf_{\epsilon \downarrow 0 } g_{\epsilon}(x_{\epsilon})   \ge g(x).  
$$
\end{definition}

Notice that Mosco convergence implies $\Gamma_s$-convergence, 
but the reverse implication is not necessarily true. 

A straightforward corollary of \cref{defn:mosco} is that both notions of convergence preserve
 convergence of minimizing sequences; more specifically for Mosco-convergence: 
\begin{proposition}\label{prop:preserve}

Suppose $x^{\ast}_{\epsilon} \rightharpoonup x^{\ast}$ and 
that $\forall \epsilon >0$, $\displaystyle x_{\epsilon}^{\ast} \in \argmin_{x\in X} g_{\epsilon}(x).$ 
If $g$ is proper and $$g_{\epsilon} \MoscoConv g,$$ 
then $ x^{\ast} \in \argmin_{x \in X} g(x)$. 
\end{proposition}
An analogous result holds when a sequence of minimizers $(x_{\epsilon}^{\ast})_{\epsilon \in \RR_+}$ 
strongly converges to $x^{\ast}$ and $(g_{\epsilon})_{\epsilon \in \RR_+}$ $\Gamma_s$-converges to $g$. 
Both modes of convergence preclude
situations like that of the simple example from earlier in this subsection, 
where the limiting point of a sequence of minimizers to 
\eqref{eq:1d_corrupted} 
is not a minimizer of  
\eqref{eq:1d_pure}, 
even when the two underlying probability distributions are arbitrarily close.

Finally, we introduce the notion of a Rockafellian associated to an optimization 
problem. 
\begin{definition}[Rockafellian]
For Banach spaces $X$ and $Y$, $\varphi: X \to \RRext$ and generic optimization 
problem $\min_{x \in X} \varphi(x)$, a bivariate function 
$\Phi: X\times Y \to \RRext$ is a Rockafellian for the problem, anchored at $\overline{y}$, 
when 
$$
\Phi(x,\overline{y}) = \varphi(x) \qquad \forall x \in X. 
$$
\end{definition}
Note that Rockafellians are not unique; for a given optimization problem, there are infinitely 
many Rockafellians associated to it. This flexibility enables one to design Rockafellians that 
are able to more easily absorb approximations than the true problem.

\section{Rockafellian relaxation and its convergence theory}
\label{sec:general_theory}

Having introduced the necessary definitions and assumptions, 
we describe in this section Rockafellian relaxation theory for 
PDECO problems of the form \eqref{eq:og_objective}. 
In particular, we consider two distinct 
scenarios---corruptions to probability densities and corruptions to 
the support of a probability distribution---and prove $\Gamma_s$-convergence
results. We also consider 
a special example of the former case, namely when the probability 
space is both discrete and finite dimensional, for which we prove
a Mosco convergence result.

   \subsection{Corruptions to continuous probability distributions} 
   \label{subsec:continous_theory}
%
%
The first type of corruption that we consider is that of a continuous 
probability distribution. 
We assume throughout this subsection that (i) $\PP$ in \eqref{eq:og_objective} is 
a probability measure on the measurable
space $(\Xi, \cA)$, (ii) there exists another sigma-finite  
measure $\mu$ on $(\Xi, \cA)$, and (iii) $\PP$ is 
absolutely continuous with respect to $\mu$. Letting  
$\rho := d\PP/d\mu$ 
denote Radon-Nikodym derivative, after a change of variables the 
objective function from \eqref{eq:og_objective} can then be written as
\begin{equation}\label{eq:objective_with_density}
\varphi(z) =  f_0(z) + \int_{\Xi} g(s(\xi, z)) \rho(\xi) \, d \mu(\xi).
\end{equation}
Taking 
$$
\Rho := \Big\{ \rho: \Xi \to \RR_+ \, \big\vert  \, \rho\in L^{\infty}(\Xi; \RR) \text{ and }  \int_{\Xi} \rho(\xi) \, d\mu(\xi) = 1 \Big\} 
$$
to denote the set of probability densities on $\Xi$,
consider
for some ``corrupted'' distribution $\rho_{\epsilon} \in \Rho$ indexed by $\epsilon >0$ 
a corresponding corrupted objective functional 
\begin{equation}\label{eq:corrupted_obj_continuous}
\varphi_{\epsilon}(z) := f_0(z) + \int_{\Xi} g(s(\xi, z)) \rho_{\epsilon}(\xi) \, d \mu(\xi) . 
\end{equation}
Additionally, for some $q \in (1,\infty)$, let 
$$
T := L^q(\Xi, \cA, \mu) = \Big\{ f: \Xi \to \RR \text{ measurable} \,\, \Big\vert \, \int_{\Xi} |f(\xi)|^q \, d\mu(\xi) < \infty \Big\}
$$
and for $q=\infty$, let $T:= L^{\infty}(\Xi,\cA,\mu)$.

Next, define a Rockafellian $\Phi: Z\times T\to \RRext$ for \eqref{eq:objective_with_density} as 
\begin{equation} \label{eq:rock_discrete_phi}
\Phi(z,t) :=
\begin{cases}
 f_0(z) + \int_{\Xi} g(s(\xi, z)) (\rho(\xi)+t(\xi)) \, d \mu(\xi) \qquad & t(\xi) = 0 \text{ a.s.} \\
\infty, & \text{else,}
\end{cases}
\end{equation}
which is of course anchored at $t(\xi) = 0$. 
Note that both here and throughout this subsection, ``a.s.'' is with respect to $\mu$.
Finally, 
define the bivariate functional $\Phi_{\epsilon}: Z \times T \to \RRext$ by setting
\begin{equation}\label{eq:mollified_rock_discrete_phi}
\Phi_{\epsilon}(z,t) := f_0(z) + \int_{\Xi} g(s(\xi, z)) (\rho_{\epsilon}(\xi)+t(\xi)) \, d \mu(\xi)
+ \frac{\theta_{\epsilon}}{q} \| t\|^q_{L^q} + \iota_{\Rho}(\rho_{\epsilon} + t).
\end{equation}
Note that $\Phi_{\epsilon}$ is parameterized by $\theta_{\epsilon} > 0$. In the special case when $q = \infty$, we instead
replace $\| t\|_{L^q}^q/q$ in \eqref{eq:mollified_rock_discrete_phi} with $\| t\|_{L^{\infty}}$. Here the
 indicator function $\iota_{\Rho}$ is defined analogously to \eqref{eq:indic_defn}. Notice that 
this is a Rockafellian for the corrupted objective functional \eqref{eq:corrupted_obj_continuous}
anchored at $t(\xi) = 0$.

We are now ready to show $\Gamma_s$-convergence of $\Phi_{\epsilon}$ to $\Phi$ whenever 
the corrupted probability density $\rho_{\epsilon}$ converges to the original, ``uncorrupted'' 
density $\rho$. 
\begin{theorem} \label{thm:density_thm}
For fixed $1\le q \le \infty$, let $(\rho_{\epsilon})_{\epsilon \in \RR_+} \subset \Rho \cap T$ 
be a sequence of probability densities that converge in $T$ to $\rho \in \Rho \cap T$. 
Additionally, let $(\theta_{\epsilon})_{\epsilon \in \RR_+} \subset \RR_+$ with $\theta_{\epsilon} \to \infty$.  
If $q \in [1,\infty)$, assume
$$
\lim_{\epsilon \downarrow 0} \theta_{\epsilon} \| \rho_{\epsilon} - \rho\|^q_{L^q} = 0,
$$
and if $q = \infty$, assume
$$
\lim_{\epsilon \downarrow 0} \theta_{\epsilon} \| \rho_{\epsilon} - \rho\|_{L^{\infty}} = 0.
$$
Finally, assume that for all $z \in Z$, $g(s(\cdot, z)) \in L^{r}(\Xi, \cA, \mu)$, where
$1/r + 1/q = 1$. 
Then, under \cref{asmpt:1} and \cref{asmpt:2}, we have 
$$
\Phi_{\epsilon} \GammaConvs \Phi
$$ 
as $\epsilon \downarrow 0$, where $\Phi$ and $\Phi_{\epsilon}$ are defined by 
\eqref{eq:rock_discrete_phi} and \eqref{eq:mollified_rock_discrete_phi}, respectively. 
\end{theorem}
\begin{remark}\label{rmk:weak_in_z}
By definition, the consistency and stability conditions in \cref{defn:mosco} in the present setting 
involve strong convergence in $Z \times T$. The proof below does rely on strong convergence in $T$; however, 
 only weak convergence in the control space $Z$ is needed. 
\end{remark}
\begin{remark}\label{rmk:theta}
The theorem states that $\Gamma_s$-convergence holds so long as $\theta_{\epsilon}$ grows
more slowly than $\| \rho_{\epsilon}-\rho\|_{L^q}^q$ vanishes. 
When employing Rockafellian relaxation in practice, one must choose an actual value for $\theta_{\epsilon}$
and, crucially, the size of the corruption $\| \rho_{\epsilon}-\rho\|_{L^q}$ is of course unknown. 
Following the numerical examples presented in \cref{subsec:numerical_examples} below, 
\cref{subsec:theta} discusses the effect of different choices for the $\theta_{\epsilon}$ parameter in
practice.  
\end{remark}
\begin{proof}[Proof of \cref{thm:density_thm}]
Here we assume that $q \in [1,\infty)$, as the case of $q=\infty$ proceeds similarly. 

To first establish the limit superior condition, suppose that 
$(z,t) \in Z \times T$. Since the condition trivially holds whenever
$t(\xi) \ne 0$ a.s., suppose $t(\xi) = 0$ a.s.
Constructing the strongly converging sequences as $z_{\epsilon} = z$
and $t_{\epsilon} = \rho - \rho_{\epsilon}$ for all $\epsilon > 0$, we have 
$$
\Phi_{\epsilon}(z_{\epsilon}, t_{\epsilon}) = f_0(z)  
+ \int_{\Xi} g(s(\xi, z)) \rho(\xi) \, d \mu(\xi)
+ \frac{\theta_{\epsilon}}{q} \| \rho - \rho_{\epsilon} \|^q_{L^q} + \underbrace{\iota_{\Rho}(\rho)}_{= 0 }
$$
and thus 
$$
\lim_{\epsilon \downarrow 0} \, \Phi_{\epsilon} (z_{\epsilon}, t_{\epsilon}) = 
f_0(z)  
+ \int_{\Xi} g(s(\xi, z)) \rho(\xi) \, d\mu(\xi)
+ \underbrace{\lim_{\epsilon \downarrow 0}  \frac{\theta_{\epsilon}}{q} \| \rho - \rho_{\epsilon} \|^q_{L^q}}_{=0}
= \Phi(z,0)
$$
as desired. 

Next, consider an arbitrary sequence $(\zeps, \teps)_{\epsilon \in \RR_+}$ that converges strongly to some $(z,t)\in Z \times T$ in the 
product topology, and first suppose that $t(\xi) = 0$ a.s.\ Since 
$(\Phi_{\epsilon}(z_{\epsilon},t_{\epsilon}))_{\epsilon\in\RR_+}$
is a sequence of extended real numbers, there necessarily exists some 
subsequence $(\Phi_{\epsilon'}(z_{\epsilon'},t_{\epsilon'}))_{\epsilon'\in\RR_+}$
such that 
\begin{equation}\label{eq:magic_start}
\lim_{\epsilon' \downarrow 0} \Phi_{\epsilon'}(z_{\epsilon'},t_{\epsilon'}) 
= \liminf_{\epsilon \downarrow 0}  \Phi_{\epsilon}(z_{\epsilon},t_{\epsilon}).
\end{equation}
Since $\rho_{\epsilon'} + t_{\epsilon'} \to \rho$ in $T$ by supposition, 
and because convergence in $T$ implies pointwise a.s.\ convergence up to a subsequence 
\cite[Exercise 5.9]{wheeden1977measure}, 
there exists a further subsequence indexed by $\epsilon''$ such that 
$\rho_{\epsilon''}(\xi) + t_{\epsilon''}(\xi) \to \rho(\xi) $
a.s.\ in $\Xi$. 

The convergence $z_{\epsilon''} \to z$ in $Z$ in turn implies both 
$\liminf_{\epsilon'' \downarrow 0} f_0(z_{\epsilon''}) \ge f_0(z)$
and 
$$
\liminf_{\epsilon'' \downarrow 0} g\big(s(\xi,z_{\epsilon''})\big) \ge  g\big(s(\xi,z)\big) \qquad \text{for } \xi\in \Xi \text{ a.s.}, 
$$
by the weak lsc property of $g$ and $f_0$ and \cref{asmpt:1}.2.

Applying Fatou's lemma (note that $g$ is bounded below by \cref{asmpt:1}.3), the above implies    
\begin{align}
\liminf_{\epsilon'' \downarrow 0}  \Phi_{\epsilon''}(z_{\epsilon''},t_{\epsilon''}) 
&\ge f_0(z) + \int_{\Xi} g\big(s(\xi,z)\big) \rho(\xi) \, d\mu(\xi) + \liminf_{\epsilon \downarrow 0 } \frac{\thetaeps}{q} \| \teps \|_{L^q}^q  \\
&\ge \Phi(z, 0);   \nonumber
\end{align}
note the final inequality comes from the trivial observation that 
$$\liminf_{\epsilon'' \downarrow 0 } \thetaeps \| \teps\|_{L^q}^q/q \ge 0.  $$
This gives the desired condition, since
$$
\liminf_{\epsilon'' \downarrow 0}  \Phi_{\epsilon''}(z_{\epsilon''},t_{\epsilon''}) 
= 	\lim_{\epsilon'' \downarrow 0}  \Phi_{\epsilon''}(z_{\epsilon''},t_{\epsilon''})
= 	\lim_{\epsilon' \downarrow 0}  \Phi_{\epsilon'}(z_{\epsilon'},t_{\epsilon'})
= 	\liminf_{\epsilon \downarrow 0}  \Phi_{\epsilon}(z_{\epsilon},t_{\epsilon})
$$
by elementary properties of sequences of real numbers. 

It remains to consider an arbitrary sequence $(\zeps, \teps)_{\epsilon \in \RR_+}$ that converges strongly 
to some $(z,t)\in Z \times T$ with $t(\xi) \ne 0$ a.s.\
In this case, the Rockafellian $\Phi(z,t) = \infty$. 
As in the above, we employ a subsequence along which $\rho_{\epsilon'}(\xi) + t_{\epsilon'}(\xi)
\to \rho(\xi) + t(\xi)$ a.s.; however, for simplicity we abuse notation and revert the subsequence
index back to $\epsilon$.
Since $\thetaeps \to \infty$ and $\lim_{\epsilon\downarrow 0} \|\teps\|_{L^q} \ne 0$ 
implies $\| \teps\|_{L^q}$ is necessarily bounded below for $\epsilon$ sufficiently small, we have 
$$
\liminf_{\epsilon\downarrow 0} \frac{\thetaeps}{q} \| \teps\|_{L^q}^q = \infty. 
$$
The desired condition
$\liminf_{\epsilon\downarrow 0} \Phi_{\epsilon}(\zeps , \teps) \ge \Phi(z,t)$
 then follows from the fact that $\iota_{\Rho}$ is bounded below and $f_0$ is proper lsc, as 
well as the estimate
\begin{align*}
\liminf_{\epsilon\downarrow 0} \int_{\Xi} g\big(s(\xi, \zeps)\big)(  \rho_{\epsilon}(\xi)+\teps(\xi) ) d\mu(\xi)
&\ge 
\int_{\Xi} g\big(s(\xi,z)\big) (\rho(\xi) + t(\xi) ) d\mu(\xi) \\
&\ge  \gamma + \int_{\Xi} g\big(s(\xi,z)\big) t(\xi)  d\mu(\xi) > -\infty, 
\end{align*}
which follows from Fatou's lemma, 
pointwise a.s.\ convergence of $\rho_{\epsilon}(\xi)+t_{\epsilon}(\xi)$,
the fact that $g$ is lsc,
and H\"{o}lder's inequality. 
\end{proof}
We remark that in the case $q = \infty$, the proof simplifies slightly; convergence in the 
$L^{\infty}$-norm topology is stronger than pointwise a.s.\ convergence, and hence working with 
subsequences to establish the limit inferior condition is not necessary.


   \subsection{Corruptions to finite, discrete probability distributions} 
   \label{subsec:discrete_theory}
   Suppose now that the underlying sample space $\Xi$ is both discrete and finite, 
so that $\PP$ in \eqref{eq:og_objective} is a finite superposition of Dirac measures. 
In this setting, we can upgrade the $\Gamma_s$-convergence from 
\cref{thm:density_thm} to the stronger notion of Mosco convergence.

For some $N \in \mathbb{N}$, redefine
\begin{equation}\label{eq:prob_vec_defn}
\Delta := \Big\{ q \in \RR^N: \, \, \sum_{i=1}^N q_i = 1  \text{ and } 0 \le q_i \le 1  \text{ for } 1 \le i \le N \Big\}
\end{equation}
be the set of probability vectors,  and let 
$\Xi = \{\xi_i\}_{i=1}^N$, where $\xi_i \in \RR^d$ for each $1 \le i \le N$ and $d \in \mathbb{N}$. 
Let 
 $\cA = 2^{\Xi}$, and for some $p \in \Delta$, let $\PP[\xi = \xi_i] = p_i$ for each $1 \le i \le N$. 

We consider 
corruptions
to the discrete probabilities $\{p_i\}_{i=1}^N$. 
This situation is relevant, for example, when 
pursuing an SAA to a quantity-of-interest
in the presence of corrupted data.

The uncorrupted objective functional \eqref{eq:og_objective} in this setting becomes
\begin{equation}\label{eq:uncorrupted_discrete_psi}
\varphi(z) := f_0(z) + \expval{ g\big( s(\xi, z)\big)} = f_0(z) + \sum_{i=1}^N p_i g\big( s(\xi_i, z)\big) .
\end{equation}
Considering now a corruption to $p$ in the form of $p_{\epsilon} \in \Delta$ (indexed by $\epsilon > 0$), the analogue to 
\eqref{eq:uncorrupted_discrete_psi} is then 
\begin{equation}\label{eq:corrupted_discrete}
\varphi_{\epsilon}(z) := f_0(z) + \sum_{i=1}^N p_{\epsilon,i} g\big( s(\xi_i, z)\big). 
\end{equation}
(here $p_{\epsilon,i}$ denotes the $i$th component of $p_{\epsilon} \in \RR^N$).

Next, define a Rockafellian 
$\Phi: Z \times \RR^N \to \RRext$
for \eqref{eq:uncorrupted_discrete_psi} as
\begin{equation}\label{eq:rock_discrete_psi}
\Phi(z,t) := f_0(z) + \sum_{i=1}^N (p_i+t_i) g\big( s(\xi_i, z)\big) + \iota_{\{0\}}(t) , 
\end{equation}
and notice that, as before, minimizing the Rockafellian $\Phi(z,t)$ is trivially equivalent to minimizing $\varphi(z)$, as 
the only feasible choice for $t$ is the zero vector. 
Finally, for some penalty parameter $\thetaeps > 0$ and $q \in [1,\infty]$, define the bivariate Rockafellian
 functional 
$\Phi_{\epsilon} : Z \times \RR^N \to \RRext$
for the corrupted objective \eqref{eq:corrupted_discrete}
by 
\begin{equation}\label{eq:mollified_rock_discrete_psi}
\Phi_{\epsilon}(z,t) := f_0(z) + \sum_{i=1}^N (p_{\epsilon,i}+t_i) g\big( s(\xi_i, z)\big) 
+ \frac{\thetaeps}{q} \| t\|_q^q +  \iota_{\Delta}(p_{\epsilon} + t), 
\end{equation}
where $\| \cdot \|_q$ denotes the $l^q$ norm on $\RR^N$. As in the previous subsection, 
in the particular case when $q = \infty$, we replace the $\| t \|_q^q/q$ term in the above 
with $\| t \|_{\infty}$.  

Next, we show that 
the sequence
 $(\Phi_{\epsilon})_{\epsilon \in \RR_+}$ Mosco converges to the Rockafellian $\Phi$ 
whenever $p_{\epsilon} \to p$ faster than the growth penalty $\thetaeps \to \infty$.

\begin{theorem}\label{thm:discrete_thm}
For fixed $1\le q \le \infty$, let
$(p_{\epsilon})_{\epsilon \in \RR_+} \subset \Delta$
be a sequence of probability vectors that converges in $l^q(\RR^N)$ to $p \in \Delta$. 
Additionally, let $(\theta_{\epsilon})_{\epsilon \in \RR_+} \subset \RR_+$ with $\theta_{\epsilon} \to \infty$. Finally, 
if $q \in [1,\infty)$, assume
$$
\lim_{\epsilon \downarrow 0} \thetaeps \| p_{\epsilon} - p\|_q^q = 0
$$ 
and if $q = \infty$, assume
$$
\lim_{\epsilon \downarrow 0} \thetaeps \| p_{\epsilon} - p\|_{\infty} =  0.
$$
Then, under \cref{asmpt:1} and \cref{asmpt:2}, 
we have 
$$
\Phi_{\epsilon} \MoscoConv \Phi  
$$
as $\epsilon \downarrow 0$, where $\Phi$ and $\Phi_{\epsilon}$ are defined by 
\eqref{eq:rock_discrete_psi} and \eqref{eq:mollified_rock_discrete_psi}, respectively. 
\end{theorem}
The proof proceeds quite similarly to that of \cref{thm:density_thm}, and hence it can be found in 
\cref{appendix:discrete_proof}. 

We remark that in the case when the sample space $\Xi$ is discrete
and countably infinite (formally, $N = \infty$), it is difficult to guarantee Mosco convergence
of $\Phi_{\epsilon}$ to $\Phi$, as weak convergence in $l^q(\NN)$ is no longer equivalent to strong convergence
(in contrast to the case of $l^q(\RR^N)$ for finite $N$). 
However, $\Gamma_s$-convergence follows as a special case of \cref{thm:density_thm}; here $\mu$ would be the
counting measure on $(\NN, 2^{\NN})$. 


   \subsection{Corruptions to the support of a probability distribution} 
   \label{subsec:support_theory}

The final type of corruption that we consider is to the support $\Xi$ of a probability distribution $\PP$. 
In this subsection we further assume that the sample space $\Xi$ is a 
separable Banach space $(\Xi, \|\cdot\|_{\Xi})$. 
First, 
define for $q \in [1,\infty)$ the space 
$$
T := L^q(\Xi; \Xi) = 
\Big\{\phi:\Xi \to \Xi \text{ strongly measurable} \,\,\Big\vert \,\, \expval{ \|\phi(\xi)\|_{\Xi}^q} < \infty \Big\}.
$$
We also can consider the case $q = \infty$, for which $T := L^{\infty}(\Xi; \Xi)$. 
For some ``corruption map'' $\etaeps \in T$
indexed by $\epsilon >0$, 
consider a corruption to the original objective functional 
$\varphi(z)$ from \eqref{eq:og_objective} in the form of 
\begin{equation}\label{eq:corrupted_objective}
\varphi_{\epsilon}(z) := f_0(z) + \expval{ g(s(\etaeps(\xi), z))}. 
\end{equation}

Next, define a 
Rockafellian functional for \eqref{eq:og_objective} by 
\begin{equation}\label{eq:pure_rock_objective}
\Phi(z,t) :=  
\begin{cases} 
f_0(z) + \expval{g\big(s(\xi+t(\xi),z)\big)}, \qquad &t(\xi) = 0 \text{ a.s.} \\
\infty &\text{else}.
\end{cases}
\end{equation}
Once again, note that minimizing the Rockafellian $\Phi(z,t)$ is trivially equivalent to 
minimizing \eqref{eq:og_objective}, as 
$t(\xi) = 0$ a.s.\ is the only feasible choice. 
Finally, for some
$\thetaeps > 0$, consider for $q \in [1,\infty)$
a Rockafellian functional $\Phi_{\epsilon}: Z\times T \to \RRext$
for the corrupted objective \eqref{eq:corrupted_objective} as 
\begin{equation}\label{eq:corrupted_rock}
\Phi_{\epsilon}(z,t) :=  
f_0(z) + \expval{g\big(s(\etaeps(\xi)+t(\xi),z)\big)} + \frac{\thetaeps}{q} \| t\|_T^q. 
\end{equation}
In the case $q = \infty$, $\| t \|_T^q/q$ in the above is replaced with $\| t \|_T$.

%

To achieve $\Gamma_s$-convergence of $\Phi_{\epsilon}$ to $\Phi$ as the size of the corruption vanishes, i.e.\ 
as $\etaeps \to I$ in $T$ (where $I(\xi) = \xi$ is the identity map), we require a stronger continuity assumption 
on the solution operator $s: \Xi \times Z \to U$ than in \cref{asmpt:1}. 

\begin{assumption}[Properties of the solution map $s=s(\xi,z)$]\label{asmpt:4}
\begin{enumerate}
\item $s(\cdot, z): \Xi \to U$ is 
$\mathcal{A}$ measurable $\forall z \in Z$.

\vskip0.5em
\item 
If both $\xi_{\epsilon} \to \xi$ in $\Xi$ and
$\zeps \rightharpoonup z$ in $Z$ as $\epsilon \downarrow 0$, 
then 
$$
s(\xi_{\epsilon}, z_{\epsilon}) \rightharpoonup s(\xi, z) \text{ in } U. 
$$ 
\end{enumerate}
\end{assumption}

A final assumption needed is that the identity map $I$ is an 
element of $T$,\footnote{We remark this is only possible if 
$\Xi$ is separable; see \cite[Example 1.1.5]{tuomas2016analysis}.}
which amounts to assuming that
the random parameter $\xi$ either has finite $q$-th moment (for $q \in [1,\infty)$) or 
is uniformly bounded in the $\Xi$ norm (for $q =\infty$). 
\begin{theorem}\label{thm:support_thm}
Fix $1 \le q \le \infty$, and suppose $(\theta_{\epsilon})_{\epsilon \in \RR_+} \subset \RR_+$
with 
$\thetaeps \to \infty$. Suppose that the identity map $I \in T$, and 
let $(\etaeps)_{\epsilon \in \RR_+} \subset T$. If $q \in [1,\infty)$ assume
$$
\lim_{\epsilon\downarrow 0} \thetaeps \| \eta_{\epsilon} - I\|_T^q = 0,
$$
and if $q = \infty$, assume
$$
\lim_{\epsilon\downarrow 0} \thetaeps \| \eta_{\epsilon} - I\|_T = 0.
$$
Then, under \cref{asmpt:2} and \cref{asmpt:4}, 
$$
\Phi_{\epsilon} \GammaConvs \Phi, 
$$
where 
$\Phi$ and $\Phi_{\epsilon}$ are 
defined by \eqref{eq:pure_rock_objective} and \eqref{eq:corrupted_rock}, respectively. 
\end{theorem}
The proof proceeds quite similarly to that of \cref{thm:density_thm}; as with the proof of 
\cref{thm:discrete_thm}, it hence can be found in 
\cref{appendix:support_proof}. 
Similar to \cref{rmk:weak_in_z}, we note that the proof 
here technically only requires weak convergence in the control space $Z$.


\section{Elliptic PDECO problems under uncertainty}
\label{sec:elliptic_pdes}
We consider in this section the optimal control of linear, elliptic partial differential 
equations (PDEs) with random coefficients in the presence of distributional ambiguity.

The goal here is to first verify that the assumptions made throughout \cref{sec:general_theory}
indeed hold for problems of this form; hence, the 
Mosco and $\Gamma_s$-convergence theorems will be applicable in this setting.
We then illustrate that these problems can be ill-conditioned
to outlier data points (similar to the simple example in \cref{subsec:motivation})
and showcase 
the ability of Rockafellian relaxation to recover the optimal controls for 
the corresponding uncorrupted problems.
\Cref{subsec:theta} and \cref{subsec:cost} conclude the
section with a brief discussion of some of the practical aspects of 
the Rockafellian relaxation, namely  
the computational cost and
the role of the 
regularization parameter $\thetaeps$ present in the objective functions 
\eqref{eq:mollified_rock_discrete_phi}, 
\eqref{eq:mollified_rock_discrete_psi}, and 
\eqref{eq:corrupted_rock}. 

\subsection{Verification of assumptions}\label{subsec:verify_assumptions}

We first introduce the model problem. 
Let $(\Xi, \cA, \PP)$ be a probability space, and
for $n \in \NN$, let $\Omega \subset \RR^n$ be a bounded, Lipschitz domain with boundary $\partial \Omega$.
Let $Z = L^2(\Omega)$ and $U = H^1_0(\Omega)$.
Given some $\alpha > 0$ and target function $u_{\star} \in L^2(\Omega)$,
we consider stochastic optimal control problems of the form 
\begin{equation}\label{eq:elliptic_stochastic_opt_cont}
\min_{z \in L^2(\Omega)} \varphi(z), \qquad \varphi(z) =  \frac12\EE\Big[  \|  s(\xi,z) - u_{\star} \|_{L^2(\Omega)}^2 \Big]
+ \frac{\alpha}{2} \| z \|_{L^2(\Omega)}^2,  
\end{equation}
where $u(\cdot, \xi):= s(\xi,z)$ 
is the solution to the elliptic PDE constraint
\begin{align}
-\nabla \cdot \big( a(x,\xi) \nabla u(x,\xi) \big) = z(x),& \qquad (x,\xi) \in \Omega\times \Xi \label{eq:2d_elliptic} \\
u(x,\xi) = 0,& \qquad (x,\xi) \in \partial\Omega\times \Xi.     \nonumber
\end{align}
Here the problem coefficient 
$a: \Omega \times \Xi \to \RR$, and 
we assume that $a(x,\cdot)$ is measurable for all $x \in \Omega$, 
and that
there exist real numbers 
$a^{\ast}, a_{\ast}  > 0$ such that 
\begin{equation}\label{eq:a_bounded_coercive}
0 < a_{\ast} \le a(x,\xi) \le a^{\ast} < \infty
\end{equation}
almost everywhere (a.e.)~in $\Omega$ and a.s.\ in $\Xi$.
\begin{definition}[Solution operator $s$]\label{defn:solution_op_elliptic}
Assume \eqref{eq:a_bounded_coercive} holds. Define the solution operator 
$s: \Xi \times L^2(\Omega) \to H^1_0(\Omega)$ to be the map that takes
some given $z\in L^2(\Omega)$ and outputs the unique solution $u(\cdot, \xi):= s(\xi,z) \in H^1_0(\Omega)$
to the variational problem: find $u(\cdot, \xi)$ such that  
\begin{equation}\label{eq:elliptic_weak_form_as}
\int_{\Omega} a(x,\xi) \nabla u(x,\xi) \cdot \nabla v(x) \, dx 
= \int_{\Omega} z(x) v(x) \, dx, \qquad \forall v \in H_0^1(\Omega). 
\end{equation}
\end{definition}
Note that existence and uniqueness of $s(\xi,z)$ follow from the Lax-Milgram theorem \cite[Corollary 8.11]{arbogast:2008}. 

We now verify \cref{asmpt:1} in the present context; in this case, it is straightforward 
to show that we actually have strong convergence of the solution operator. 
\begin{proposition}\label{prop:elliptic_1}                                       
Suppose $(\zeps)_{\epsilon \in \RR_+} \subset L^2(\Omega)$ and 
$\zeps \rightharpoonup z$ in $L^2(\Omega)$.
Then 
\begin{equation}\label{eq:desired_weak_conv}
s(\xi,\zeps)\to s(\xi, z) \text{ in } H^1_0(\Omega) 
\end{equation}
a.s.\ in $\Xi$. 
\end{proposition}
The proof follows from standard arguments, but, for completeness, it is included in \cref{appendix:verify_proofs}. 

To next verify \cref{asmpt:1}.1, note that 
the Lax-Milgram theorem also guarantees there exists a unique solution $u \in L^2(\Xi; H^1_0(\Omega))$
such that
\begin{equation}\label{eq:elliptic_weak_form_exp}
\EE\Big[\int_{\Omega} a(x,\xi) \nabla u(x,\xi) \cdot \nabla v(x,\xi) \, dx \Big] = \EE\Big[\int_{\Omega} z(x) v(x,\xi) \, dx\Big]
\end{equation}
for all $v \in L^2(\Xi; H_0^1(\Omega))$, which in particular guarantees that the solution map is measurable.  
Note that the solution to \eqref{eq:elliptic_weak_form_exp} coincides with the 
solution to \eqref{eq:elliptic_weak_form_as}.

Using the notation of \eqref{eq:og_objective}, 
the control regularizer and quantity-of-interest map in \eqref{eq:elliptic_stochastic_opt_cont} 
are 
$$
f_0(z) := \frac{\alpha}{2} \| z\|^2_{L^2(\Omega)} \qquad \text{and} \qquad  g(u) := \frac12  \| u - u_{\star} \|^2_{L^2(\Omega)},
$$
respectively, which clearly fulfill the requirements of \cref{asmpt:2}.1 and \cref{asmpt:2}.3. 
Hence, \cref{asmpt:2} holds for 
the stochastic elliptic optimal control problem. 
Ergo, \cref{thm:discrete_thm} follows in the context described in \cref{subsec:discrete_theory}.
By the Lax-Milgram theorem, we know that for any $z \in L^2(\Omega)$, $g(s(\cdot, z)) \in L^{\infty}(\Xi, \cA, \PP)$. 
Hence, assuming that $\mu(\Xi) < \infty$
(in the context of \cref{subsec:continous_theory}), \cref{thm:density_thm} will hold.

It remains to verify the 
additional hypotheses on the solution 
map $s(\xi, z)$ in \cref{asmpt:4} of \cref{subsec:support_theory}; for this 
an additional continuity property on the 
random coefficient $a(x,\xi)$ is needed, which we describe now.  
As in \cref{prop:elliptic_1}, we can actually show strong convergence. 

\begin{proposition}\label{prop:elliptic_2}                                               
Suppose that the sample space is a separable Banach space $(\Xi, \| \cdot \|_{\Xi})$.
Suppose $(\zeps)_{\epsilon \in \RR_+} \subset L^2(\Omega)$ and 
$\zeps \rightharpoonup z$ in $L^2(\Omega)$.
Suppose $(\xi_{\epsilon})_{\epsilon \in \RR_+} \subset \Xi$ and 
$\xi_{\epsilon} \to \xi$ in $\Xi$, and assume \eqref{eq:a_bounded_coercive}
holds at both $\xi$ and $\xi_{\epsilon}$ for all $\epsilon\in\RR_+$. 
If $a(x,\cdot): \Xi \to \RR$ is sequentially continuous for almost every $x$ in $\Omega$, then
\begin{equation}\label{eq:desired_weak_conv_2}
 s(\xi_{\epsilon}, \zeps) \to s(\xi,z) \text{ in } H^1_0(\Omega). 
\end{equation}
\end{proposition}
As before, the proof can be found in 
\cref{appendix:verify_proofs}.
Notice that
continuity of $a(x,\xi)$ is only needed in $\xi$; 
discontinuities in the $x$ variable are permitted, 
so long as \eqref{eq:a_bounded_coercive} holds. 
Such continuity holds, for example, when $a(x,\xi)$ is given by a Kosambi-Karhunen-Lo\'{e}ve
(KKL) approximation\footnote{The polymath Damodar Dharmananda Kosambi 
invented the expansion in 1943 \cite{kosambi1943}, preceding both 
Karhunen (1946, \cite{karhunen1946spektraltheorie})
and Lo\'{e}ve (1948, \cite{loeve1948}).} \cite[Theorem 2.3]{martinez2018optimal}
of a log-normal random field 
$$
a(x,\xi)  = e^{m(x) + \sigma(x) \sum_{k=1}^d \sqrt{\lambda_k} b_k(x) \xi_k  }. 
$$ 
  For each $1 \le k \le d$, $(\lambda_k, b_k(x))$ are eigenvalue-eigenfunction 
pairs of the integral operator defined by a covariance kernel $C(x,x')$ associated
to some Gaussian field. The functions $m(x)$ and $\sigma(x)$ are the mean and variance that result from 
the covariance function, and in turn $e^{m(x)}$ are $e^{\sigma(x)}$ are the resulting geometric
mean and geometric variance of the log-normal random field. 
Finally, each $\xi_k$ are independent and
normally distributed with zero mean and unit variance. See \cite[Remark 2.4]{martinez2018optimal}
for a complete description. 


   \subsection{Numerical examples}
   \label{subsec:numerical_examples}
   Having verified the assumptions necessary for the convergence theorems described
in \cref{sec:general_theory} to hold in the case of stochastic optimal control 
problems constrained by linear, elliptic PDEs,  
we now illustrate (i) the sensitivity of this class of problems to perturbations 
of the underlying probability densities and sample space and (ii) the ability of 
Rockafellian relaxation to recover an optimal control $z^{\ast}$ to an \emph{uncorrupted}
problem in the presence of data corruption. We describe three examples. 

In all examples, we consider the problem of minimizing the objective function 
\eqref{eq:elliptic_stochastic_opt_cont}, but in the first two examples, 
the objective is constrained by the following one-dimensional elliptic boundary value problem (BVP)
posed on the domain $\Omega = (0,1)$:
\begin{align} 
-\frac{d}{dx} \Big( a(x,\xi) \frac{d}{dx} u(x,\xi) \Big) = z(x)&, \qquad (x,\xi) \in (0,1)\times \Xi \label{eq:1d_bvp_constraint} \\
u(0,\xi) = u(1,\xi) = 0&.  \nonumber
\end{align}
The control regularization parameter in \eqref{eq:elliptic_stochastic_opt_cont} 
is set as $\alpha = 10^{-4}$. 
The state equation \eqref{eq:1d_bvp_constraint} and its
corresponding adjoint are both discretized
with a standard second-order 
finite difference method on a uniform grid 
with a mesh spacing $h = 1/256$, and the solutions to the resulting tridiagonal linear systems are 
computed with SciPy's banded direct solver \cite{harris2020array}. 
The $L^2$ norms $\| \cdot \|_{L^2(0,1)}$ are computed
with the composite trapezoidal rule.

In the third example, we take $\alpha = 10^{-5}$ and constrain
 the objective function 
\eqref{eq:elliptic_stochastic_opt_cont}
by the two-dimensional elliptic BVP \eqref{eq:2d_elliptic}
posed on the domain $\Omega = B_1(0)$, i.e.\ the unit ball in $\RR^2$
centered at the origin. The state and adjoint equations are solved with a 
finite element method with $Q1$ elements implemented in \texttt{deal.II} 
\cite{arndt2022deal}; the total number of degrees of freedom is 5185.

The specific form
of the diffusion coefficient $a(x,\xi)$,  
the sample space $\Xi$, and the probability distribution associated to the 
the random parameter $\xi$ will vary across each example, as we detail below.


\vskip0.5em
\textit{Example 1:}  
Consider first setting $a(x,\xi) = \xi$, so that 
\eqref{eq:1d_bvp_constraint} becomes a scaled Poisson problem. 
In the uncorrupted case, let $\xi$ be the discrete random variable
defined by 
\begin{equation}\label{eq:uncorrupted_xi_1}
\PP[\xi = 2] = 1, 
\end{equation}
which simply results in a deterministic optimal control problem. In the corrupted case, 
define the random variable $\xi_{\epsilon}$ by 
\begin{equation}\label{eq:corrupted_xi_1}
p_{\epsilon,1} := \PP[\xi_{\epsilon} = 0.2] = \epsilon, \qquad p_{\epsilon,2}:=  \PP[\xi_{\epsilon} = 2] = 1-\epsilon, 
\end{equation}
where $0 < \epsilon \ll 1$.  Let the target function $u_{\star}(x) = \sin(\pi x)$. 
The minimization problem in each case is solved using a standard 
gradient descent method 
with backtracking line search based on the Armijo condition. 
The gradient descent is terminated when the inner product of the control and 
the gradient of the objective functional is less than $\tau = 10^{-4}$. 
The initial guess 
is set as $z(x) = 1$. 

This problem is analogous to the simple one-dimensional stochastic program described in 
\cref{subsec:motivation}, in the sense that a minimizer $z^{\ast}_{\rm corrupted}$ to the 
corrupted problem is quite different from a minimizer $z^{\ast}_{\rm true}$ to the 
uncorrupted one, even when $\epsilon$ is small. Indeed, even for a 
0.5\% corruption (i.e.\ for $\epsilon = 0.005$), the pointwise difference between 
$z^{\ast}_{\rm true}$ and $z^{\ast}_{\rm corrupted}$ is $\bigoh{1}$ throughout most 
of the domain. 

\begin{figure}[h]       
    \centering
(a)
 \includegraphics[width=0.45\textwidth]{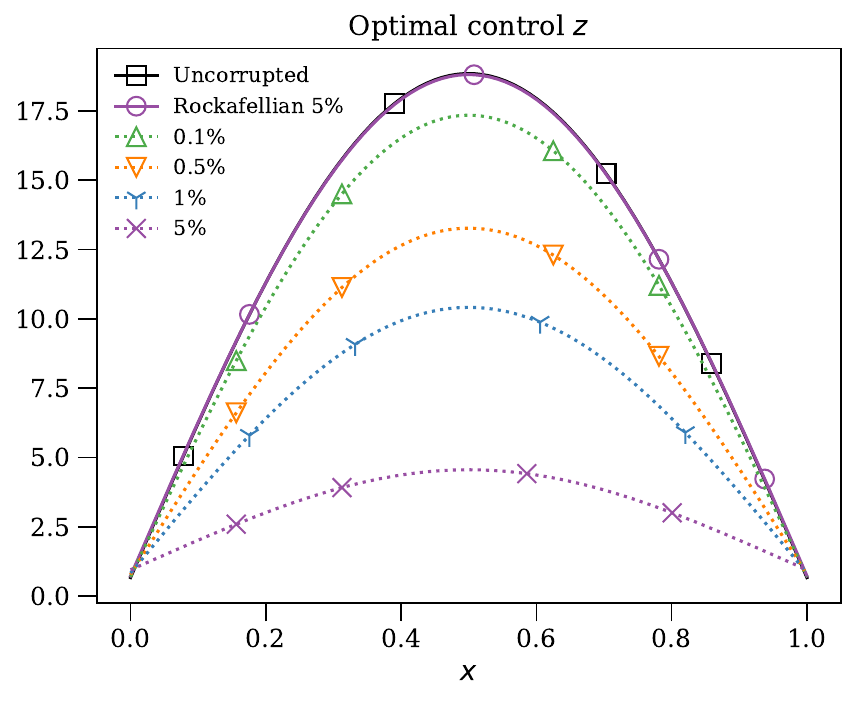}
(b)
 \includegraphics[width=0.45\textwidth]{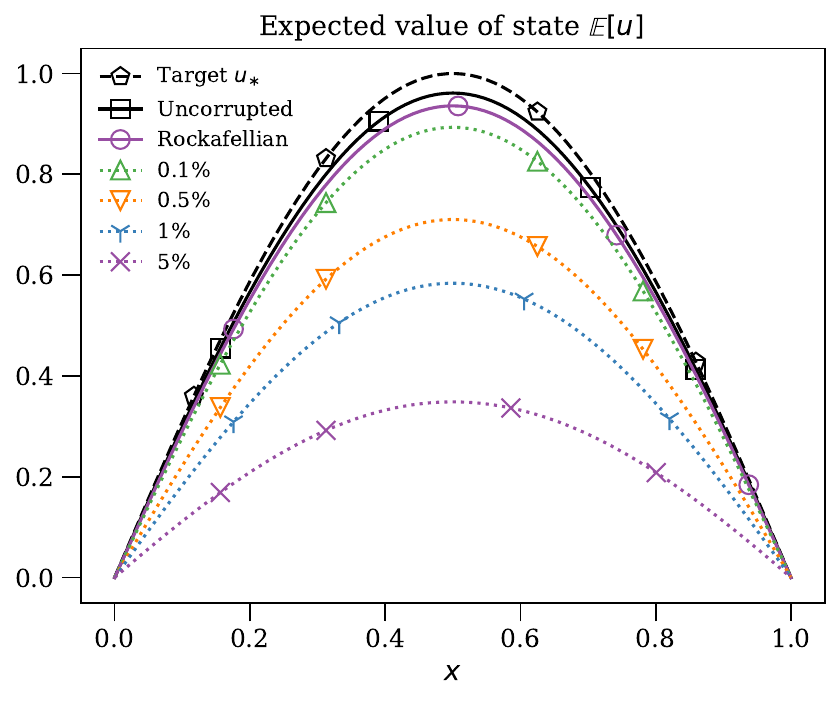}
\caption{Example 1: 
(a) Optimal controls for the true, uncorrupted problem, the corrupted
problem at varying 
corruption levels (dotted lines), and the Rockafellian relaxation at 5\% corruption. (b) The 
expected value of the corresponding solutions to the state 
equation \eqref{eq:1d_bvp_constraint}
in each case. Corruption levels are defined by \eqref{eq:corrupted_amount_defn_case1}.
} 
\label{fig:plots_case1}
\end{figure}

\cref{fig:plots_case1}(a) shows $z^{\ast}_{\rm true}$ and $z^{\ast}_{\rm corrupted}$ 
for varying corruption levels
\begin{equation}\label{eq:corrupted_amount_defn_case1}
\text{\% corruption} := 100\cdot \epsilon,
\end{equation}
while \cref{fig:plots_case1}(b) shows the expected value of the corresponding solutions 
to the state equation \eqref{eq:1d_bvp_constraint}. In contrast to the simple example 
in \cref{subsec:motivation}, the optimal control for the corrupted problem 
 appears to converge 
to that of the uncorrupted one as $\epsilon \downarrow 0$ (as seen in \cref{fig:plots_case1}); 
however, the convergence 
is slow, as there are nontrivial errors even at 0.1\% corruption. 

Consider now the bivariate Rockafellian  
$\Phi_{\epsilon}: L^2(\Omega) \times \RR^2 \to \RRext$,
\begin{align}
\Phi_{\epsilon}(z,t):= \frac12 
\sum_{i=1}^2 (p_{\epsilon,i} &+ t_i) \| s(\xi_{\epsilon,i}, \, z) - u_{\star} \|_{L^2(0,1)}^2 \label{eq:rock_obj_ex1}
+ \frac{\alpha}{2} \| z\|_{L^2(0,1)}^2 \\
&+ \frac{\theta}{2} \| t\|^2_2 + \iota_{\Delta} (p_{\epsilon}+t),  \nonumber 
\end{align}
where $\theta > 0$, $\xi_{\epsilon,1} = 0.2$, $\xi_{\epsilon,2} = 2$, and the indicator function $\iota_{\Delta}$
constrains $p_{\epsilon} + t$ to be a probability 
vector in $\RR^2$ (see \eqref{eq:prob_vec_defn}), which is equivalent to the constraint that 
$t_1 + t_2 = 0$ 
and $-p_{\epsilon,i} \le t_i \le 1-p_{\epsilon,i}$, $i\in\{1,2\}$. 

In practice we enforce the equality constraint 
by eliminating $t_1$ via 
$t_1 = -t_2$. 
The Rockafellian $\Phi_{\epsilon}(z,t_2)$ is then optimized 
with the projected gradient descent method with backtracking line search based on the Armijo condition; 
here the projection enforces the bound constraints 
$-p_{\epsilon,2} \le t_2 \le 1-p_{\epsilon,2}$.
The initial guesses are $z(x) = 1$ (as before) and $t_2=0$. 

\cref{fig:plots_case1}(a) shows the resulting optimal control $z^{\ast}_{\rm Rock}$ 
for $\epsilon = 0.05$ (5\% corruption) and $\theta = 1$. 
By optimizing the relaxed Rockafellian $\Phi_{\epsilon}$, we recover
a minimizer $z^{\ast}_{\rm true}$ to the uncorrupted problem, as desired; 
the absolute error 
$$
\| z^{\ast}_{\rm true} - z^{\ast}_{\rm Rock} \|_{L^{2}(0,1)} \le 
\| z^{\ast}_{\rm true} - z^{\ast}_{\rm Rock} \|_{L^{\infty}(0,1)} = 5.55\cdot 10^{-2}.
$$ 
Similar results (using the same $\theta$ value) are obtained for the other corruption levels (0.1\%, 0.5\%, or 1\%).


\vskip0.5em
\textit{Example 2:}                                                      
For $d \in \NN$ and $\sigma \in \RR_+$, consider next setting the coefficient
$$
a(x,\xi) = 
e^{\sigma\sum_{k=1}^{d}  \sqrt{\lambda_k} \sin(x/\sqrt{\lambda_k}) \xi_k  }, \qquad \lambda_k = \frac{4}{(2k-1)^2\pi^2}, 
$$ 
where each $\xi_k$ are independent and normally distributed with zero mean and unit variance; 
here the argument in the exponential is a truncated KKL expansion for a rescaled Brownian motion 
on the interval $\Omega=(0,1)$. 

The expectation value in \eqref{eq:elliptic_stochastic_opt_cont} is approximated with an SAA
using $N$ samples; the objective function then becomes
\begin{equation}\label{eq:obj_fn_saa}
\varphi(z) =  \frac12 \sum_{i=1}^{N} p_i \|  s(\xi^{(i)},z) - u_{\star} \|_{L^2(0,1)}^2 
+ \frac{\alpha}{2} \| z \|_{L^2(0,1)}^2,
\end{equation}
where $p_i = 1/N$ for all $i$ and the samples $\xi^{(i)} \sim \mathcal{N}(0,I_{d})$ are all independent (here $I_{d}$ is the 
$d\times d$ identity matrix). Define the target function as $u_{\star}(x) = 1$. 

In all cases, we take $d = 50$ terms in the truncated KKL expansion and $N=1000$ SAA
samples. 
In the uncorrupted case, we take $\sigma = 0.4$, while in the corrupted 
cases, we select the first $M \in \NN$ samples ($M < N$) and increase their variance
through the map $\xi \mapsto 10\,\xi$, which amounts to setting $\sigma = 4$ for those samples. 
Optimization is done with the SciPy \cite{2020SciPy-NMeth} 
implementation of the BFGS algorithm; the algorithm terminates 
whenever the norm of gradient of the objective function is less than $10^{-5}$, 
which is the default setting. 
In all cases, the initial guess $z(x) = 1$. 

\begin{figure}[h]       
    \centering
(a)
 \includegraphics[width=0.45\textwidth]{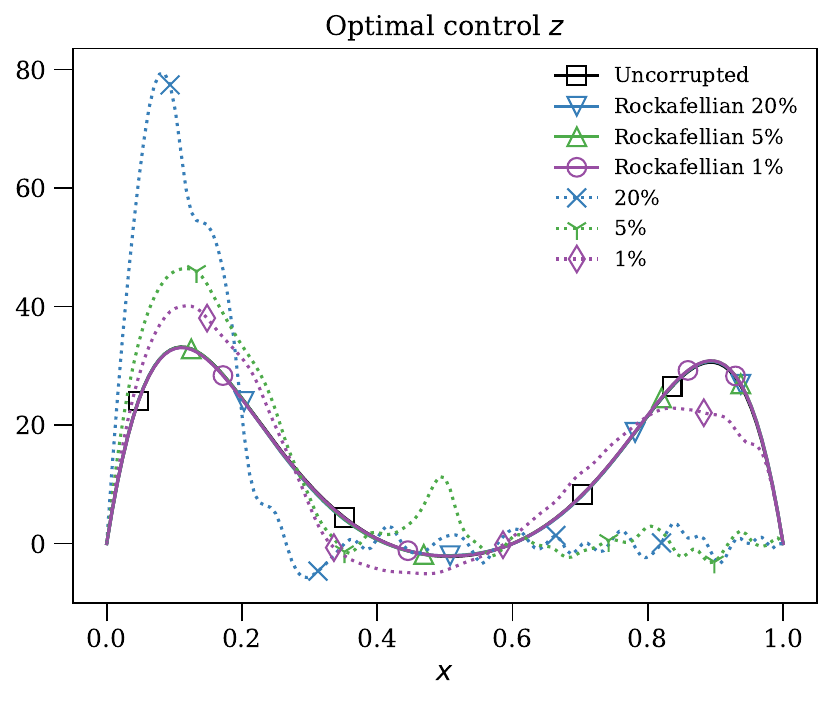}
(b)
 \includegraphics[width=0.45\textwidth]{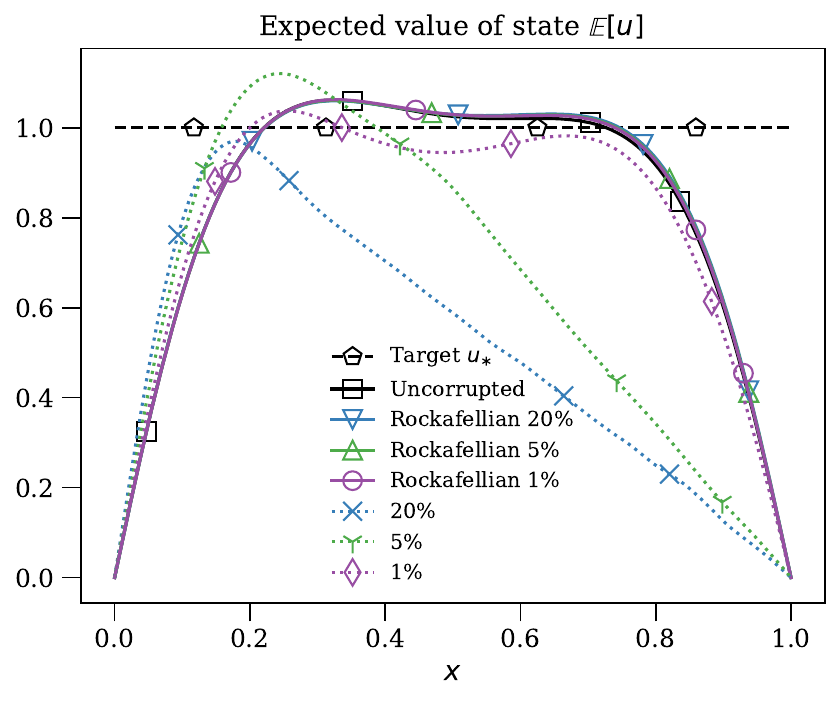}
\caption{Example 2:  
(a) Optimal controls for the true, uncorrupted problem, as well as for the corrupted
problem (dotted lines) and the corresponding Rockafellian relaxations at varying 
corruption levels. (b) The 
expected value of the corresponding solutions to the state 
equation \eqref{eq:1d_bvp_constraint}
in each case. Corruption levels are defined by \eqref{eq:corrupted_amount_defn_case4}.
} 
\label{fig:plots_case4}
\end{figure}

With the corruption level defined as
\begin{equation}\label{eq:corrupted_amount_defn_case4}
\text{\% corruption} := 100 \cdot  M/N,
\end{equation}
\cref{fig:plots_case4}(a) shows that 
even at just 1\% corruption, there are severe differences
between $z^{\ast}_{\rm corrupted}$ and $z^{\ast}_{\rm true}$, 
i.e.\ minimizers of
\eqref{eq:obj_fn_saa} with and without corruption, respectively. 
These differences grow with increasing corruption levels. 
\cref{fig:plots_case4}(b) shows the corresponding expected values of the corrupted state variables 
are also quite different from the uncorrupted one.

To recover $z_{\rm true}^{\ast}$ in the presence of corrupted data samples, define the bivariate 
Rockafellian $\Phi: L^2(0,1)\times \RR^N \to \RRext$: 
\begin{equation}\label{eq:obj_fn_saa_rock}
\Phi(z,t) :=  \frac12 \sum_{i=1}^{N} (p_i+t_i) \|  s(\xi^{(i)},z) - u_{\star} \|_{L^2(0,1)}^2 
+ \frac{\alpha}{2} \| z \|_{L^2(0,1)}^2 + \theta \| t\|_1 + \iota_{\Delta} (p+t).
\end{equation}
As the role of the perturbation variable here is to identify (and remove) outlier 
sample points, $t$ is measured in the $l_1$ norm for its well known 
sparsity-promoting property; 
this choice is inspired by the statistical learning experiments described in 
\cite[Section 6]{royset2024rockafellian}, 
where corrupted labels in the contexts of computer vision and text analytics 
were properly removed. 

The experiments in \cite{royset2024rockafellian} also inspired the optimization 
method for \eqref{eq:obj_fn_saa_rock};  
an ADI heuristic is adopted 
in which we first fix $t=0$ and compute $z^{\ast} \in \argmin \Phi(z,0)$ using the BFGS method. 
The result is then used to compute the solution to the linear 
program $t^{\ast} \in \argmin \Phi(z^{\ast},t)$ 
using SciPy's implementation of the simplex method. This process is repeated until the absolute $l_1$ 
distance 
between successive $t^{\ast}$ values is smaller than some $\tau$; in particular we take $\tau=10^{-5}$, 
consistent with the BFGS stopping criteria quoted in Example 1 above.

The linear program solve in the ADI approach 
naturally allows the $p+t \in \Delta$ constraint to be satisfied
exactly; in particular this constraint imposes the 
pointwise bounds $-p_i \le t_i \le 1-p_i$ for each $1 \le i \le N$.
In practice, to prevent the linear program solver from ``greedily'' identifying a handful of 
samples $\xi^{(i)}$ at which 
$$
\| s(\xi^{(i)}, z) - u_{\star}\|_{L^2(0,1)}^2
$$
is the smallest and subsequently deleting all of the other samples (even those that are ``clean'', i.e. 
uncorrupted), 
we additionally impose the more stringent pointwise bounds $-p_i \le t_i \le p_i$
for each $i$. 

\cref{fig:plots_case4}(a) shows the optimal controls $z^{\ast}_{\rm Rock}$ for \eqref{eq:obj_fn_saa_rock} 
with $\theta = 5\cdot 10^{-2}$ at 1\%, 2\%, and 10\% corruption levels; 
 the true optimal control $z^{\ast}_{\rm true}$ in each case is indeed recovered, as desired. 
\cref{fig:plots_case4}(b) shows that the corresponding expected values of the state variables in the Rockafellian
case are much more accurate than the corrupted counterparts.  
\cref{table:case_4_accuracy} shows the relative $L^2$ errors in the optimal controls, as well as the ratio 
of errors for the corrupted and Rockafellian cases, as defined by 
\begin{equation} \label{eq:rel_l2_err_defn}
\qquad E_{\rm rel}(z) := 
\frac{
\| z  - z_{\rm true}^{\ast} \|_{L^2(0,1)} }
{\| z_{\rm true}^{\ast} \|_{L^2(0,1)}},
\qquad 
\mathcal{E}_{\rm ratio} := \frac{E_{\rm rel}(z^{\ast}_{\rm corrupted})}{E_{\rm rel}(z^{\ast}_{\rm Rock})} . 
\end{equation} 
The errors for $z_{\rm Rock}^{\ast}$ are all one to two orders of magnitude 
smaller than the corresponding ones in the corrupted cases, i.e.\ $\mathcal{E}_{\rm ratio}$
is large. 
Although not shown in \cref{fig:plots_case4}, the table also 
shows that Rockafellian relaxation can even recover the optimal controls 
at a 40\% corruption level.\footnote{We note that this, and even 20\% and 10\% 
corruption levels, may be unrealistically large in many engineering contexts.} 
In all cases, the Rockafellian perturbation variable $t$
identifies and removes at least 70\% of the corrupted data points, 
while less than 3\% of the clean, uncorrupted points
are incorrectly removed. 
In general these percentages will change as a function 
of $\theta$, which we discuss below in \cref{subsec:theta}.

\begin{table}                                   
  \begin{center}
  \def~{\hphantom{0}}
    \begin{tabular}{ c | c | c | c | c  }
\hline
Corruption level & $E_{\rm rel}(z^{\ast}_{\rm Rock})$ & $\mathcal{E}_{\rm ratio}$ & Corrupted deleted & Clean deleted \\
\hline
  1\%            & $6.84\cdot 10^{-3}$ & 37.4  &  7/10=70\% &  26/990=2.62\%                     \\
  5\%            & $7.72\cdot 10^{-3}$ & 97.8 &  40/50=80\% &   23/950=2.42\%                     \\  
  20\%            & $8.96\cdot 10^{-3}$ & 115.4 &  171/200=85.5\% & 16/800=2.00\%                     \\
  40\%           & $2.20\cdot 10^{-2}$ & 47.3 &  336/400=84\% & 11/600=1.83\%                     \\
    \end{tabular}
    \caption{
Example 2 with $\theta = 5\cdot 10^{-2}$: relative $L^2$ errors $E_{\rm rel}$  
between the Rockafellian and true optimal controls, as well as 
the ratio of $L^2$ errors $\mathcal{E}_{\rm ratio}$ for the corrupted and 
Rockafellian optimal controls; see \eqref{eq:rel_l2_err_defn}. Also included
are the fraction of corrupted and 
clean sample points correctly and incorrectly removed by the perturbation 
variable $t$, respectively. 
Corruption levels defined by \eqref{eq:corrupted_amount_defn_case4}. 
}
  \label{table:case_4_accuracy}
  \end{center}
\end{table}


\vskip0.5em
\textit{Example 3:}                                                       
For the final example, we consider corruptions to the 
support $\Xi$ to a probability distribution $\PP$, as 
in \cref{subsec:support_theory}. The two-dimensional PDE constraint is 
given by  
\eqref{eq:2d_elliptic} on the domain $\Omega = B_1(0)$. We take 
\begin{equation}\label{eq:acoeff_case5}
a(x,\xi) = \frac{1}{ \xi + 3 \sin(10\pi \|x\|_2) }, 
\end{equation}
and in the uncorrupted case, 
$\xi = 3.5$ almost surely.
For the corrupted cases, 
$\xi$ is uniformly distributed on the interval 
$[3.5-\delta,3.5+\delta]$, where $\delta\in (0,0.5)$.
Notice that that the contrast ratio 
\begin{equation}\label{eq:contrast}
\sup_{x \in \Omega} a(x,\xi)/\inf_{x \in \Omega} a(x,\xi)
\end{equation}
of the oscillations
is large when $\xi$ takes on values close to $3$. 

Since the Dirac measure is not suitable for the 
theory developed in \cref{subsec:support_theory}, 
for conceptual purposes
one could instead take the uncorrupted random variable $\xi$ to be
uniformly distributed on the interval $[3.5-\nu,3.5+\nu] \subset \Xi := \RR$, 
where $0 < \nu \ll 1$. The corruption map $\eta : \Xi\to\Xi$ would
then be 
$
\eta(\xi) = (\delta/\nu) \xi + 3.5(1-\delta/\nu),
$
which of course tends
to the identity as $\delta \searrow \nu$. 

The expectation value in \eqref{eq:elliptic_stochastic_opt_cont} 
is discretized with standard Gaussian 
quadrature with $N = 8$ points, and optimization in both the uncorrupted
and corrupted cases is done with 
a L-BFGS optimizer with history size $m = 7$ \cite[Chapter 7.2]{nocedal1999numerical}
 and backtracking line search based on the Wolfe criterion.
The target function $u_{\star}(x) = 1$, the initial guess $z(x) = 1$, and 
$\alpha = 10^{-5}$.

\begin{figure}[h]
   \begin{center}
       \includegraphics[width=0.30\textwidth]{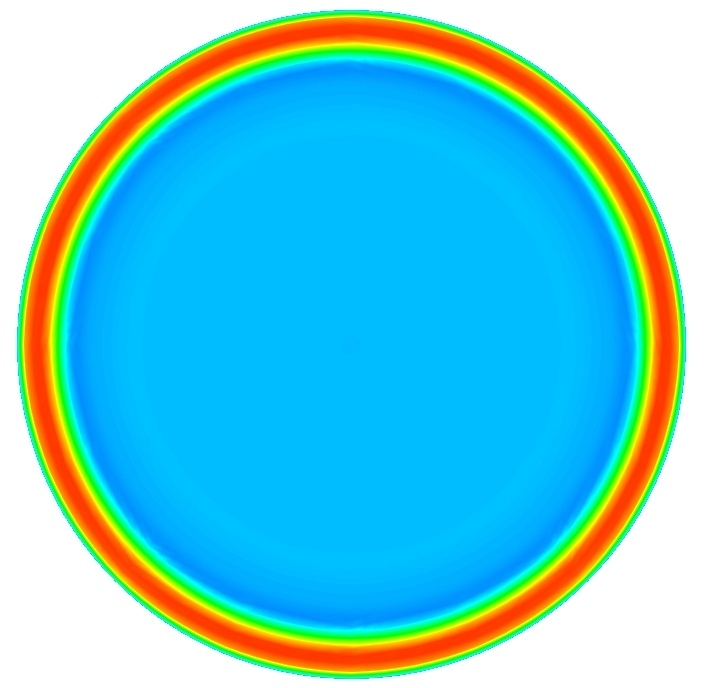}
       \includegraphics[width=0.30\textwidth]{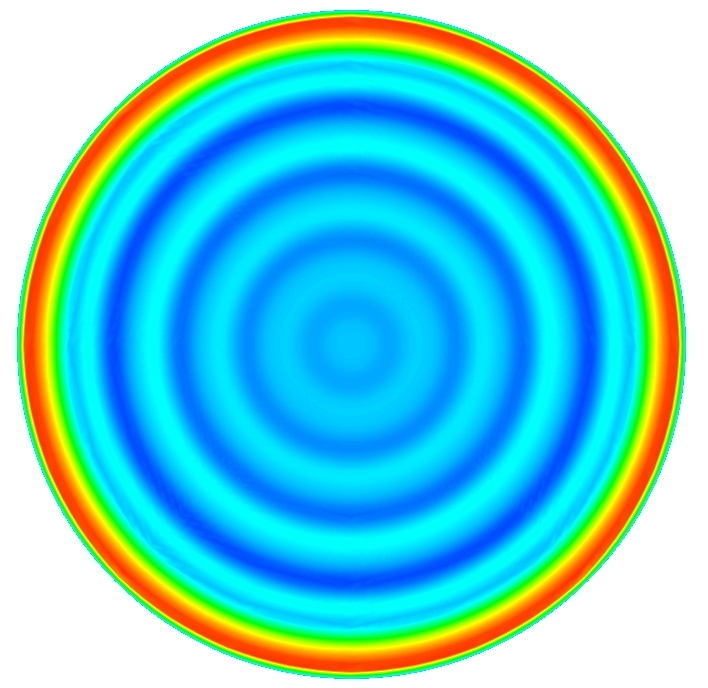}
       \includegraphics[width=0.30\textwidth]{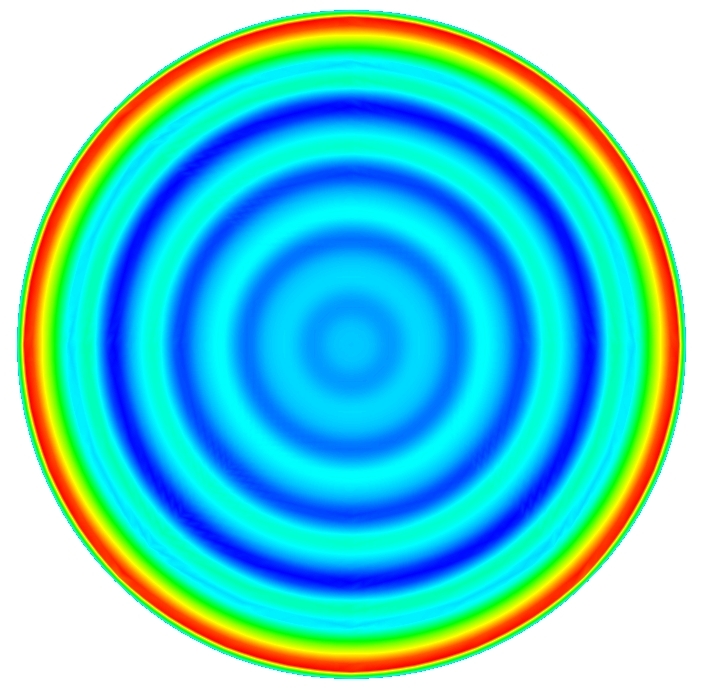}
       \includegraphics[width=0.06\textwidth]{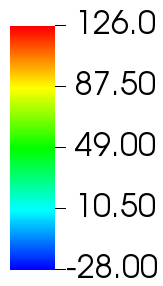} 
   \end{center}
\caption{
Example 3: optimal control for the true, uncorrupted problem (left), as well as 
for the corrupted problems with $\delta = 0.3$ (middle) and $\delta = 0.4$ (right).
}
\label{fig:plots_case5}
\end{figure}

\cref{fig:plots_case5} shows that for $\delta = 0.3$ and $\delta = 0.4$, the oscillations
in the coefficient $a(x,\xi)$ lead to oscillatory optimal controls. 
The amplitude of the oscillations increases as $\delta$ gets larger. 

To recover an optimal control $z_{\rm true}^{\ast}$ for the uncorrupted problem 
in the presence of increased uncertainty in $\xi$, 
we optimize the following bivariate Rockafellian $\Phi: L^2(\Omega)\times L^2(\Xi;\Xi) \to \RR$:
\begin{align}\label{eq:rock_ex_3}
\Phi(z,t) := \frac{1}{2} \int_{\Xi_{\delta}} \| s(\tilde\xi &+t(\tilde\xi),z) - u_{\star} \|_{L^2(\Omega)}^2 \rho(\tilde\xi)\, d\tilde \xi 
+ \frac{\alpha}{2} \| z\|^2_{L^2(\Omega)} \\
&+ \frac{\theta}{2} \int_{\Xi_{\delta}} |t(\tilde\xi)|^2 \rho(\tilde\xi)\, d\tilde\xi,   \nonumber
\end{align}
where $\theta > 0$, $\Xi_{\delta} = [3.5-\delta,3.5+\delta]$, $\rho(\tilde \xi) = \mathds{1}_{\Xi_{\delta}}(\tilde \xi)/(2\delta)$, 
and $\mathds{1}_S$ denotes the 
characteristic function on the set $S$; here we use $\tilde \xi$ to 
denote the corrupted random variable. 
We additionally enforce the bound 
constraints on the Rockafellian variable $t(\tilde \xi)$: 
$$
3.5-\delta \le \tilde \xi + t(\tilde \xi) \le 3.5 + \delta, 
$$
so that the superposition $\tilde \xi + t(\tilde \xi)$ does not take 
values outside of the support of $\rho(\tilde \xi)$.

As in Example 2, the Rockafellian \eqref{eq:rock_ex_3} is 
optimized with an ADI heuristic. With the initial guess of 
$t(\tilde \xi) = 0$, we first compute $z^{\ast} \in \argmin \Phi(z,0)$ using the
L-BFGS method and subsequently compute $t^{\ast} \in \argmin \Phi(z^{\ast},t)$ using 
 projected gradient descent with backtracking line search based on the Armijo
criterion. This process is repeated until the absolute distance in $L^2(\Xi; \Xi)$
between successive $t^{\ast}$ values is smaller than $\tau = 10^{-2}$. 
Note that 
the Fr\'{e}chet derivative of the objective $\Phi$ with respect 
to $t$ is given by
\begin{equation}\label{eq:dPhidt}
\frac{\delta \Phi}{\delta t}(\tilde \xi) = \theta \,  t(\tilde \xi)  
- \int_{\Omega} \frac{\partial}{\partial t} a\big(x,\tilde\xi + t(\tilde \xi) \big) \nabla u(x,\tilde \xi) \cdot \nabla p(x,\tilde \xi) \, d\tilde \xi,
\end{equation}
where $p(x,\tilde \xi)$ is the solution to the adjoint equation 
\begin{align*} 
-\nabla \Big( a\big(x,\tilde\xi+t(\tilde \xi)\big) \nabla p(x,\tilde\xi) \Big) 
= u(x,\tilde \xi) - u_{\star}(x)&, \qquad (x,\tilde\xi) \in \Omega \times \Xi_{\delta}  \\
 p(x,\tilde\xi) = 0&, \qquad (x, \tilde\xi) \in \partial \Omega \times \Xi_{\delta}.  \nonumber
\end{align*}
This calculation assumes, of course, differentiability of the coefficient $a(x,\tilde \xi)$  with respect 
to $\tilde \xi$,  which is true for \eqref{eq:acoeff_case5} on the domain $\Xi_{\delta}$. 

\begin{figure}
   \begin{center}
       \includegraphics[width=0.30\textwidth]{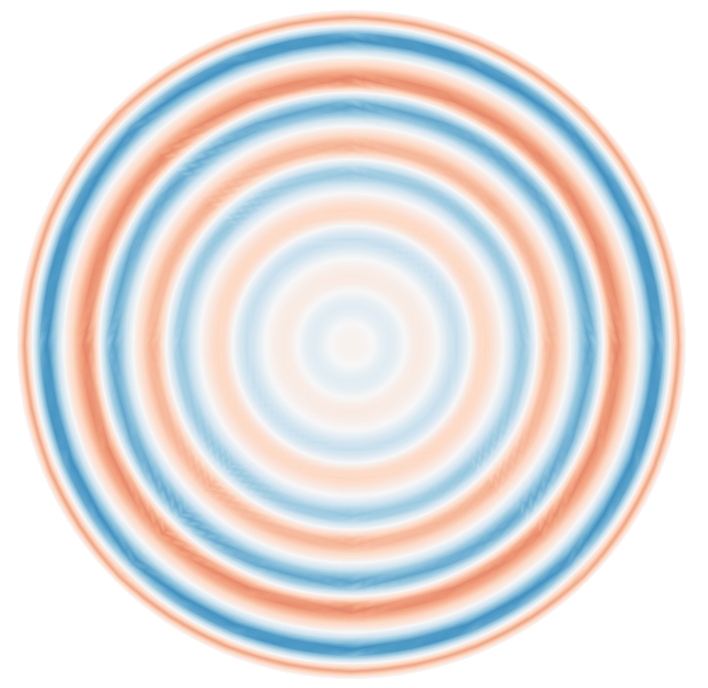}
       \includegraphics[width=0.30\textwidth]{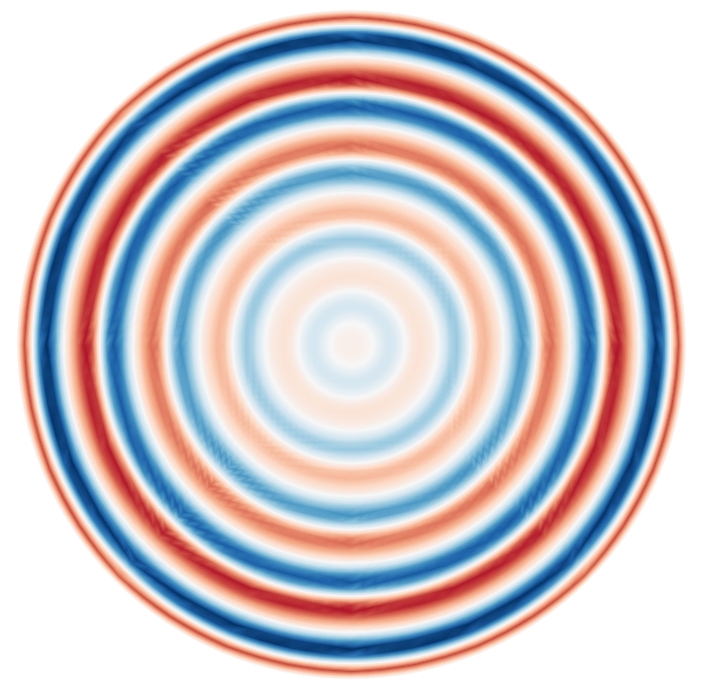}
       \includegraphics[width=0.30\textwidth]{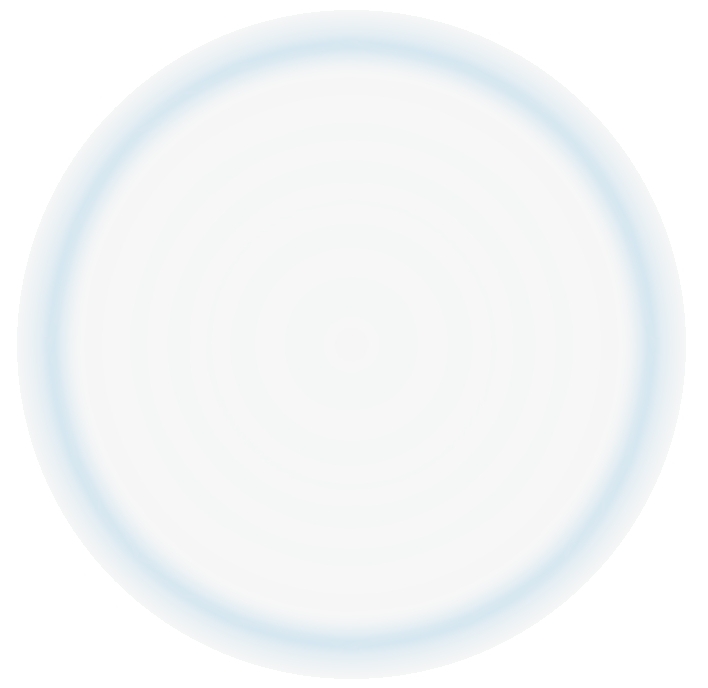}
       \includegraphics[width=0.06\textwidth]{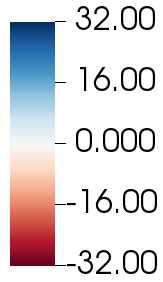} 
   \end{center}
\caption{
Example 3:
pointwise error between an uncorrupted optimal control
and corrupted optimal controls 
for $\delta = 0.3$ and $\delta = 0.4$ (left and middle, respectively), 
as well as the pointwise error for a Rockafellian optimal control for $\delta = 0.4$ (right). 
}
\label{fig:error_case5}
\end{figure}

\cref{fig:error_case5} shows that the pointwise error between an optimal control $z^{\ast}_{\rm Rock}$ for 
\eqref{eq:rock_ex_3} at $\theta = 10^{-1}$ and $\delta = 0.4$ is considerably smaller than the pointwise error in 
a corrupted optimal control $z^{\ast}_{\rm corrupted}$. In particular, the oscillations seen in \cref{fig:plots_case5}
are greatly suppressed. The $L^{\infty}$ error in $z^{\ast}_{\rm Rock}$ is approximately six times smaller than 
that of $z^{\ast}_{\rm corrupted}$.  
The results are similar when $\delta = 0.3$ (not shown).

\cref{table:case_5_accuracy}
lists the relative $L^2$ errors (defined by \eqref{eq:rel_l2_err_defn}) in the 
Rockafellian optimal controls, which are
a factor of more than six (respectively, eight) times smaller than the relative errors for 
the corrupted optimal controls 
at $\delta= 0.3$ (resp., $\delta = 0.4$).

\begin{table}                                   
  \begin{center}
  \def~{\hphantom{0}}
    \begin{tabular}{ c | c | c |  c  }
\hline
$\delta$ & $E_{\rm rel}(z^{\ast}_{\rm Rock})$ & $\mathcal{E}_{\rm ratio}$ &  $\cV_{\rm ratio}$ \\
\hline
  0.3            & $3.04\cdot 10^{-2}$ & 6.31  &  $1.959\cdot 10^2$                    \\  
  0.4            & $3.58\cdot 10^{-2}$ & 8.72 &  $1.624\cdot 10^2$                   \\
    \end{tabular}
    \caption{
Example 3 with $\theta = 10^{-1}$: relative $L^2$ errors $E_{\rm rel}$  
between the Rockafellian and true optimal controls, as well as 
the ratio of $L^2$ errors $\mathcal{E}_{\rm ratio}$ for the corrupted and 
Rockafellian optimal controls; see \eqref{eq:rel_l2_err_defn}. 
Also included
is the ratio of variances in the $L^2$ norm 
of the state variables for the corrupted and Rockafellian controls, as defined by 
\eqref{eq:variances}. 
}
  \label{table:case_5_accuracy}
  \end{center}
\end{table}

Owing to larger 
contrast ratios of $a(x,\tilde\xi)$ at values of $\tilde \xi$ close to 3 (see \eqref{eq:contrast}), the 
variance in the $L^2$ norm of the corrupted state variables 
increases with increasing $\delta$. In contrast, 
the corresponding variance for the Rockafellian state variable is considerably smaller. 
This is quantified in \cref{table:case_5_accuracy}, which shows the variance ratio
\begin{equation}\label{eq:variances}
\cV_{\rm ratio}:= \frac{\text{Var}(u_{\rm corrupted})}{\text{Var}(u_{\rm Rock})},
\end{equation}
where $\text{Var}(u):= \EE[(\|u(\cdot,\xi)\|_{L^2}-\EE[\|u(\cdot,\xi)\|_{L^2}])^2]$. 
This reduction in the variance is explained by the extremely low spread in the 
values of $\tilde \xi + t^{\ast}_{\rm Rock}(\tilde\xi)$. 
Indeed, the expected values of this 
random variable range between approximately $3.57$ and $3.59$ for $\delta = 0.3$ and $\delta = 0.4$, respectively, 
while the standard deviation in both cases is smaller than $1.6\cdot 10^{-2}$; 
essentially, 
the objective function \eqref{eq:rock_ex_3} in this case 
is minimized whenever
the Rockafellian variable $t(\tilde \xi)$ alters the 
stochastic problem \eqref{eq:rock_ex_3} to be an (approximately) deterministic one. 


   \subsection{Impact of $\theta$ parameter}
   \label{subsec:theta}
   As mentioned in both \cref{rmk:theta} and 
the preceding subsection, minimizers $z^{\ast}$
and $t^{\ast}$ of Rockafellian objective 
functionals will depend on the value of the regularization parameter
$\thetaeps$.\footnote{In this subsection we abuse notation and use $\theta$  
and $\thetaeps$ interchangeably.}
\cref{thm:density_thm} stipulates that a
Rockafellian will $\Gamma$-converge to the corresponding
uncorrupted objective functional when both $\thetaeps \to \infty$ and 
$\thetaeps \| \rho_{\epsilon}-\rho\|_T^q \to 0$ (where $T$ is an
$L^q$ norm) as $\epsilon\downarrow0$; in other words, $\thetaeps$ should grow 
as the size of the corruption vanishes, but not too quickly.  
Analogous conditions are needed in \cref{thm:discrete_thm} 
and \cref{thm:support_thm}. 

These conditions are consistent with intuition; 
 $\theta$ values that are too large will 
prevent the perturbation 
variable $t$ from meaningfully altering any corrupted problem data. 
If $\theta$ is sufficiently small, however, there is the possibility 
that $t$ ``greedily'' alters the problem data too much in 
an effort to reduce the objective functional. In the context of 
corrupted empirical samples, for example (as considered in 
\cref{subsec:discrete_theory}), this corresponds to 
improperly removing ``clean'', uncorrupted samples. 

\begin{table}                                   
  \begin{center}
  \def~{\hphantom{0}}
    \begin{tabular}{ c | c | c | c | c }
\hline
 $\theta$ & $E_{\rm rel}(z^{\ast}_{\rm Rock})$ & $\mathcal{E}_{\rm ratio}$ & Corrupted deleted & Clean deleted \\
\hline
  $5\cdot 10^{-3}$            & $4.92\cdot 10^{-2}$ & 5.62 &  19/20=95\% &   416/980=42.4\%  \\
  $5\cdot 10^{-2}$            & $7.23\cdot 10^{-3}$ & 38.3 &  14/20=70\% &   25/980=2.55\%   \\  
  $5\cdot 10^{-1}$            & $1.43\cdot 10^{-2}$ & 19.4 &  3/20=15\% &   0/980=0.00\%     \\  
    \end{tabular}
\caption{
Results from Example~2 at various $\theta$ values and 2\% corruption 
(as defined by \eqref{eq:corrupted_amount_defn_case4}). $E_{\rm rel}$ and $\cE_{\rm ratio}$ 
are defined by \eqref{eq:rel_l2_err_defn}. 
}
  \label{table:case_4_theta_effect}
  \end{center}
\end{table}

This intuition is backed up by numerical experiments at varying 
$\theta$ values. As an illustration, consider
again Example~2 in \cref{subsec:numerical_examples} at 2\% corruption. 
\cref{table:case_4_theta_effect} shows that while the relative $L^2$ 
error is $\bigoh{10^{-2}}$ or smaller for $\theta$ at three different
orders of magnitude, the error is smallest when most (70\%) of the 
corrupted samples are deleted and few (less than 3\%) of the clean
sampled are removed. 

For Rockafellian relaxations on the support of a probability distribution 
(as considered in \cref{subsec:support_theory}), the value of $\theta$ 
affects the variance 
reduction property observed in Example~3 of \cref{subsec:numerical_examples}. 
Recall that the variance in 
$\tilde\xi + t_{\rm Rock}^{\ast}(\tilde\xi)$ (denoting
the superposition of the corrupted random variable 
and an optimal perturbation variable $t$)
is quite small at $\theta=10^{-1}$; for $\theta=10^{-2}$, 
the variance is even smaller, while for $\theta=1$ it is larger. 
At smaller $\theta$ values, then, 
the variance in the state variable 
$u_{\rm Rock}(x,\xi)$ is lower; the opposite is true at larger $\theta$
values. This is quantified by the ratio $\cV_{\rm ratio}$ (defined
in \eqref{eq:variances}) in \cref{table:case_5_theta_effect} 
for the case when $\delta =0.4$. 
For all values of $\theta$,
the relative $L^2$ error between the Rockafellian optimal 
control and the true, uncorrupted one is $\bigoh{10^{-2}}$. 

\begin{table}                                   
  \begin{center}
  \def~{\hphantom{0}}
    \begin{tabular}{ c | c | c |  c  }
\hline
$\theta$ & $E_{\rm rel}(z^{\ast}_{\rm Rock})$ & $\mathcal{E}_{\rm ratio}$ &  $\cV_{\rm ratio}$ \\
\hline
  $10^{-2}$    & $4.27\cdot 10^{-2}$ & 7.30  &  $1.427\cdot 10^{4}$      \\
  $10^{-1}$    & $3.58\cdot 10^{-2}$ & 8.72  &  $1.624\cdot 10^{2}$      \\  
  $1$          & $7.43\cdot 10^{-2}$ & 4.20  &  $4.075\cdot 10^{0}$    \\
    \end{tabular}
    \caption{
Results from Example~3 at various $\theta$ values and $\delta = 0.4$. 
$E_{\rm rel}$ and $\cE_{\rm ratio}$ are defined
by \eqref{eq:rel_l2_err_defn}, and $\cV_{\rm ratio}$ is defined by 
\eqref{eq:variances}. 
}
  \label{table:case_5_theta_effect}
  \end{center}
\end{table}


   \subsection{Computational cost}
   \label{subsec:cost}
   In general, the computational cost of optimizing 
a bivariate Rockafellian is larger than 
optimizing its (potentially corrupted) 
single-variable counterpart; 
the enhanced resiliency to data corruption
afforded by Rockafellian relaxation does not come for free. 

For example, the total number of gradient descent iterations
in Example~1 taken to minimize the Rockafellian  
\eqref{eq:rock_obj_ex1} is 3187, compared to only 446 iterations 
for the corrupted problem. Note, however, that 
computing the gradient of \eqref{eq:rock_obj_ex1} with respect 
to $t$ is essentially free, since the relatively expensive 
terms  
$
\| s(\xi_{\epsilon,i},z) - u_{\star}\|^2_{L^2(0,1)}
$ 
which appear in this gradient are already needed 
for an objective function evaluation. 
Hence, the additional computational cost comes essentially
from the need for additional evaluations of the objective functional 
\eqref{eq:rock_obj_ex1} and its gradient with respect to the control $z$. 
This is true for the alternating direction (ADI) heuristics employed in 
Examples 2 and 3 as well.

Indeed, the heuristic for Example~2 entails alternating between a BFGS solve 
for the control $z$ and a linear program solve for the perturbation 
variable $t$ with $\bigoh{1000}$ 
unknowns. The cost of the latter will of course increase with the 
number of discrete samples, but for large-scale PDECO problems this
will likely remain (significantly) smaller than the cost of the former.
 
The ADI heuristic for Example~3 involves a standard
projected gradient descent (with line search) for $t$, 
and 
in practice we observe this to converge
quite rapidly---always in fewer 
than ten iterations. Note that constructing the gradient
with respect to $t$ only requires calculating
an integral over the spatial domain $\Omega$ for each 
stochastic collocation point used in the discretization 
of the probability sample space; see \eqref{eq:dPhidt}.

\begin{table}                                   
  \begin{center}
  \def~{\hphantom{0}}
    \begin{tabular}{ c | c | c | c | c }
\hline
 $\theta$ &  $N_{\rm iters}$, 2\% 
          &  $N_{\rm iters}$, 20\% 
          & $N_{\rm evals}$, 2\%
          & $N_{\rm evals}$, 20\%     \\
\hline
  $5\cdot 10^{-3}$       &  537   &  1640   & 2024   & 6251         \\
  $5\cdot 10^{-2}$       &  514   &  1852   & 1953   & 7065         \\  
  $5\cdot 10^{-1}$       &  351   &  2697   & 1333   & 10303        \\  
    \end{tabular}
\caption{
Example 2: total number of BFGS iterations, denoted $N_{\rm iter}$, 
as well as the total number of evalations 
of the objective functional \eqref{eq:obj_fn_saa_rock} and its 
gradient with respect to the control $z$, denoted
$N_{\rm evals}$, at various $\theta$ values and two
corruption levels (as defined by \eqref{eq:corrupted_amount_defn_case4}).  
}
  \label{table:case_4_cost}
  \end{center}
\end{table}

Since the additional cost for Rockafellian relaxation comes primarily 
from the additional work to optimize over the control variable $z$, we
quantify how much is needed for representative cases of Examples 2 and 3
in \cref{table:case_4_cost} and 
\cref{table:case_5_cost}, respectively.

For Example~2, 
optimizing the corrupted objective functional 
\eqref{eq:obj_fn_saa} 
at 2\% corruption without Rockafellian relaxation using
the BFGS method took 189 total iterations, 
which includes 719 evaluations of both the 
objective functional and its gradient with respect to the control $z$.
\cref{table:case_4_cost} shows that 
optimizing the Rockafellian \eqref{eq:obj_fn_saa_rock} at this level 
of corruption is about two to three times more expensive, depending
on $\theta$. 
At 20\% corruption, 412 BFGS iterations and 1576 objective
functional and gradient evaluations
were required to optimize \eqref{eq:obj_fn_saa}, while 
the cost to optimize the Rockafellian ranged between approximately 
 four to six and 
a half times 
more expensive.

For Example~3, the total number of L-BFGS iterations needed to
optimize the corrupted objective functional \eqref{eq:elliptic_stochastic_opt_cont}
is 299 and 308 for $\delta = 0.3$ and $\delta = 0.4$, respectively. 
The number of iterations needed to optimize the corresponding Rockafellian 
\eqref{eq:rock_ex_3} in $z$ is somewhere between two to three times greater. 
There is little variation as a function of $\theta$; 
see \cref{table:case_5_cost} for the precise numbers.

\begin{table}                                   
  \begin{center}
  \def~{\hphantom{0}}
    \begin{tabular}{ c | c | c   }
\hline
 $\theta$ &  $N_{\rm iters}$, $\delta = 0.3$ 
          &  $N_{\rm iters}$, $\delta = 0.4$   \\
\hline
  $10^{-2}$       &  695   &  807   \\ 
  $10^{-1}$       &  706   &  777   \\ 
  $1$             &  723   &  834   \\ 
    \end{tabular}
\caption{
Example 3: total number of L-BFGS iterations (denoted $N_{\rm iter}$)
at various $\theta$ values and two
corruption levels.
}
  \label{table:case_5_cost}
  \end{center}
\end{table}

\section{Conclusions}
\label{sec:conclusions}
We introduce a framework based on Rockafellian relaxation for 
general stochastic PDE constrained optimization 
problems that is robust 
to meta-uncertainty.
Theoretical $\Gamma$-convergence results are shown, 
and 
numerical examples of elliptic PDECO problems under uncertainty
illustrate the framework's utility 
for outlier detection and removal and for variance reduction.  
These properties depend on the regularization parameter $\theta$
inherent to Rockafellian objective functionals. 
As discussed in \cref{subsec:theta}, 
in general, $\theta$ should be neither too large nor too small relative
to the size of the data corruption in the problem; an improper choice could
have potentially unfavorable properties. In practice, of 
course, the corruption level is unknown, and if possible, decision makers should 
compute for a range of $\theta$ values. 
This is consistent with the view that, on the whole, 
optimization technology is more a tool for 
identifying possibilities than an ``oracle'' for producing a definite answer 
\cite[Preface]{primer}. 

There are numerous potential avenues for future research. 
Firstly, the framework can be employed in 
a much broader range of contexts than considered here. Two 
possible examples are data assimilation problems and full waveform inversion
problems, because of the potential for corrupted empirical measurements
exists in those contexts.

It may also be reasonable to consider different forms of the 
bivariate, Rockafellian functional. The work here considered 
$L^p$ norms for the perturbation variable $t$, but other choices, 
for example the Kullback-Leibler divergence,
are possible and may have their own advantages. 

Finally, there is an opportunity in future research to develop theory for how 
to best optimize the bivariate objective functionals that arise in Rockafellian relaxation.
While the optimization techniques described \cref{subsec:numerical_examples} were successful 
in the present setting, there is room to develop more efficient methods and to prove 
convergence guarantees. 
Along these lines, an inexact trust region framework seems promising for optimally balancing computational
cost and efficiency 
\cite{ARConn_NIMGould_PhLToint_2000a,MHeinkenschloss_LNVicente_2001a,kouri2013,kouri2014}.

%
%
%
%
%
%
%
%
%

\section*{Acknowledgments}
The authors thank Rohit Khandelwal for insightful discussions, as well as 
 Mike Novack for helpful discussions on $\Gamma$-convergence.   
Additionally, 
the authors greatly appreciate the referees for their careful reviews
and many helpful suggestions that significantly improved the manuscript. 

\appendix
   \section{Proof of \cref{thm:discrete_thm}} 
   \label{appendix:discrete_proof}
   \begin{proof}
The proof for the cases $q = \infty$ and $q \in [1,\infty)$ proceed identically; 
here we assume the latter. Let $T$ denote $l^q(\RR^N)$, and note because $T$ is
finite dimensional, weak and strong convergence are equivalent.

We first
establish the limit superior condition in \cref{defn:mosco}. Suppose 
that $(z,t) \in Z \times T$; since $\Phi(z,t) = \infty$ for $t \ne 0$
(in which case the condition holds trivially), suppose $t = 0$. 
We construct the
strongly converging sequences as
$z_{\epsilon} = z$ and $t_{\epsilon} = p - p_{\epsilon}$ for all $\epsilon > 0$.
This yields
$$
\Phi_{\epsilon}(z_{\epsilon}, t_{\epsilon} ) = f_0(z) + \sum_{i=1}^{N} p_i g\big( s(\xi_i,z)\big) 
+ \frac{\thetaeps}{q} \| p - p_{\epsilon}\|_q^q + \underbrace{\iota_{\Delta}(p)}_{= 0 }, 
$$
which implies 
$$
\lim_{\epsilon \downarrow 0} \, \Phi_{\epsilon}(z_{\epsilon}, t_{\epsilon} ) = 
f_0(z) + \sum_{i=1}^{N} p_i g\big( s(\xi_i,z)\big) 
+ \underbrace{\lim_{\epsilon \downarrow 0} \frac{\thetaeps}{q} \| p - p_{\epsilon}\|_q^q}_{=0} 
= \Phi(z,0)
$$
by assumption, as desired. 

Next, we consider an arbitrary sequence $(z_{\epsilon}, t_{\epsilon})_{\epsilon}$ 
converging weakly to some $(z,t)$ in the product topology on $Z \times T$. Since
$T$ is a finite-dimensional normed space, note that weak and strong convergence are equivalent. 
Since $f_0$ is proper lsc and  
the indicator $\iota_{\Delta}$ is also
lsc (as $\Delta$ is a closed, convex set in $\RR^N$), we have
\begin{align*}
\liminf_{\epsilon \downarrow 0} \Phi_{\epsilon}(z_{\epsilon}, t_{\epsilon}) 
\ge f_0(z) 
+ \liminf_{\epsilon \downarrow 0} \sum_{i=1}^N (p_{\epsilon,i} + t_{\epsilon,i}) &g\big( s(\xi_i, z_{\epsilon})\big) \\
+ \liminf_{\epsilon \downarrow 0} \frac{\thetaeps}{q} \| t_{\epsilon}\|_q^q
&+ \iota_{\Delta} (p + t).
\end{align*}
Since $g$ too is weakly sequentially lsc and bounded below, and since and $p_{\epsilon} + t_{\epsilon} \to p + t$, 
$$
\liminf_{\epsilon \downarrow 0} \sum_{i=1}^N (p_{\epsilon,i} + t_{\epsilon,i}) g\big( s(\xi_i, z_{\epsilon})\big) 
\ge
\sum_{i=1}^N  (p_{i} + t_{i})  g\big(s(\xi_i, z)\big)
$$
which then gives  
$$
\liminf_{\epsilon \downarrow 0} \Phi_{\epsilon}(z_{\epsilon}, t_{\epsilon}) 
\ge f_0(z) + \sum_{i=1}^N (p_i + t_i) g\big(s(\xi_i, z)\big)  
+ \liminf_{\epsilon \downarrow 0} \frac{\thetaeps}{q} \| t_{\epsilon}\|_q^q
+ \iota_{\Delta} (p + t).
$$
If $t = 0$, then the right-hand side of this inequality 
is trivially greater than or equal to $\Phi(z,0)$. If $t\ne 0$, then the Rockafellian $\Phi(z,t) = \infty$, 
while the sequence $\| t_{\epsilon}\|_q$ is necessarily bounded below by some positive 
constant for $\epsilon$ sufficiently small. Since $\thetaeps \to \infty$, 
this gives 
$
\liminf_{\epsilon \downarrow 0 } \Phi_{\epsilon} (z_{\epsilon},t_{\epsilon}) = \Phi(z,t) = \infty,
$
as desired. 
\end{proof}

   \section{Proof of \cref{thm:support_thm}} 
   \label{appendix:support_proof}
   \begin{proof}
We assume that $q \in [1,\infty)$, as the case of $q=\infty$ is nearly the same. 

Establishing first the limit superior condition, suppose that $(z,t) \in Z \times T$. 
Since the condition trivially holds for $t(\xi) \ne 0$ a.s., suppose otherwise. 
Constructing the sequences as $z_{\epsilon} = z$ and $t_{\epsilon} = \etaeps - I$
for all $\epsilon > 0$, we have
$$
\Phi_{\epsilon}(\zeps, \teps) = f_0(z) + 
\expval{g\big(s(\xi,z)\big)} + \frac{\thetaeps}{q} \| \etaeps - I\|_T^q.
$$
Taking the limit of both sides gives 
$
\lim_{\epsilon \downarrow 0} \Phi_{\epsilon}(\zeps,\teps) = \Phi(z,0), 
$
as desired. 

Next, consider an arbitrary sequence $(\zeps, \teps)_{\epsilon \in \RR_+}$ that converges strongly to some $(z,t)\in Z \times T$ in the 
product topology, and first suppose that $t(\xi) = 0$ a.s.\  
As detailed in the proof of \cref{thm:density_thm}, it suffices to establish the limit inferior condition along a subsequence for which 
\begin{equation}\label{eq:pntwise_sub}
\eta_{\epsilon'}(\xi) + t_{\epsilon'}(\xi) \to \xi \qquad \text{a.s. in } \Xi.
\end{equation}
For ease of exposition, we abuse notation by reverting the subsequence index back to $\epsilon$ and 
proceed to work with the subsequence.\footnote{In the case $q = \infty$, the subsequence argument is not necessary, 
as convergence in the $L^{\infty}$ norm implies pointwise convergence a.s.} 

The weak convergence $\zeps \rightharpoonup z$ and \eqref{eq:pntwise_sub} imply by \cref{asmpt:4} that 
$s(\eta_{\epsilon}(\xi) + t_{\epsilon}(\xi), \zeps) \rightharpoonup s(\xi, z) 
$ in $U$, 
while the lsc property of $g$ implies 
$
\liminf_{\epsilon \downarrow 0} g\big( s(\eta_{\epsilon}(\xi) + t_{\epsilon}(\xi), \zeps) \big) \ge g\big( s(\xi, z)\big).
$
Since $g$ is coercive and $f_0$ is lsc, the result 
$
\liminf_{\epsilon \downarrow 0} \Phi_{\epsilon}(\zeps, \teps) \ge \Phi(z, 0) 
$
follows from Fatou's lemma. 

Finally, consider an arbitrary sequence $(\zeps, \teps)_{\epsilon \in \RR_+}$ that converges to some $(z,t)$, with $t(\xi) \ne 0$ a.s.\ 
Because
$\Phi(z,t) = \infty$, $g$ and $f_0$ are bounded below and proper lsc, respectively, the result 
$
\liminf_{\epsilon \downarrow 0} \Phi_{\epsilon}(\zeps, \teps) = \Phi(z, t)
$
follows from the fact that 
$$
\liminf_{\epsilon \downarrow 0} \frac{\thetaeps}{q} \| \teps\|_T^q = \infty, 
$$
as $\| \teps \|_T$ is necessarily bounded below and 
 $\thetaeps \to \infty$. 
\end{proof}

   \section{Proofs of propositions from \cref{subsec:verify_assumptions}} 
   \label{appendix:verify_proofs}
   
\begin{proof}[Proof of \cref{prop:elliptic_1}]
For any $\xi \in \Xi$ such that \eqref{eq:a_bounded_coercive} holds, 
the Lax-Milgram theorem implies that the solution to \eqref{eq:elliptic_weak_form_as} 
satisfies 
$$
\| s(\xi,z) \|_{H^1_0(\Omega)} \lesssim  \| z \|_{L^2(\Omega)} .
$$
Thus, the linear operator $A(\xi): H^{1}_0(\Omega) \to H^{-1}(\Omega)$ 
defined by the left-hand side of \eqref{eq:elliptic_weak_form_as}
is continuous and has a bounded inverse. Additionally, the embedding operator
$i:  L^2(\Omega) \to H^{-1}(\Omega)$ is compact \cite[Theorem 7.20]{arbogast:2008}.  
Hence, the composition $A(\xi)^{-1}\circ i$ maps the weakly convergent sequence
$(\zeps)_{\epsilon\in\RR_+}$ to a strongly convergent sequence
$(s(\xi,\zeps))_{\epsilon\in\RR_+}$. 
\end{proof}

\begin{proof}[Proof of \cref{prop:elliptic_2}]
Recall from the statement of the Proposition that we assume
$(\xi_{\epsilon})_{\epsilon\in\RR_{+}}$ is a sequence that converges to some $\xi \in \Xi$
and that \eqref{eq:a_bounded_coercive} is assumed to hold at both $\xi$ and $\xi_{\epsilon}$ 
for all $\epsilon\in\RR_{+}$.
Define now $w_{\epsilon}:\Omega \times \Xi \to \RR$ by $w_{\epsilon}(\cdot,\xi):= s(\xi,\zeps)$. By \cref{prop:elliptic_1}, we have
$\| s(\xi,\zeps) - s(\xi,z)\|_{H^1_0(\Omega)} \to 0 $ as $\epsilon \downarrow 0$. 
Using the definitions of $w_{\epsilon}(\cdot,\xi)$ and $u_{\epsilon}(\cdot, \xi_{\epsilon}):= s(\xi_{\epsilon},\zeps)$,  
we have
$$
\int_{\Omega}\big[ a(x,\xi_{\epsilon})\nabla u_{\epsilon}(x,\xi_{\epsilon}) - 
a(x,\xi)\nabla w_{\epsilon}(x,\xi)\big] 
\cdot \nabla v(x) \, dx 
= 0, \qquad \forall v \in H^1_0(\Omega), 
$$
which implies 
$$
\int_{\Omega} a(x,\xi)\nabla[u_{\epsilon}(x,\xi_{\epsilon})-w_{\epsilon}(x,\xi)]\cdot \nabla v(x) \, dx
=\int_{\Omega} [a(x,\xi) - a(x,\xi_{\epsilon})] \nabla u_{\epsilon}(x,\xi_{\epsilon})\cdot \nabla v(x) \, dx
$$
for all $v \in H^1_0(\Omega)$. Testing with $v = u_{\epsilon}(\cdot,\xi_{\epsilon})-w_{\epsilon}(\cdot,\xi)$
and using coercivity of $a$, the Cauchy-Schwartz inequality, and uniform boundedness 
of $u_{\epsilon}(\cdot,\xi_{\epsilon})$
in $H^1_0(\Omega)$ yields
 $$
\| u_{\epsilon}(\cdot,\xi_{\epsilon}) - w_{\epsilon}(\cdot,\xi) \|_{H^1_0(\Omega)} \lesssim \esssup_{x \in \Omega}|a(x,\xi)-a(x,\xi_{\epsilon})| ,
$$
which vanishes in the limit $\epsilon \downarrow 0$ by the assumed continuity of $a$. 
After applying the triangle inequality,
$$
\| s(\xi,z) -  s(\xi_{\epsilon},\zeps)\|_{H^1_0(\Omega)} \le 
 \| s(\xi,z) -  s(\xi,\zeps)\|_{H^1_0(\Omega)} + 
\| s(\xi,\zeps) -  s(\xi_{\epsilon},\zeps)\|_{H^1_0(\Omega)}
$$
the desired result follows in the limit $\epsilon \downarrow 0$.
\end{proof}


\bibliographystyle{siamplain}
\bibliography{references}
\end{document}